\documentclass[12pt,a4paper]{amsart}
\usepackage{mathrsfs}
\usepackage{amssymb}
\usepackage{amsxtra}
\usepackage{cleveref,comment}
\usepackage{graphics,color,graphicx}
\usepackage[T1]{fontenc}
\usepackage[utf8]{inputenc}
\usepackage{tikz}
\usepackage{accents}
\usepackage{multicol}
\usepackage{caption}

\newcommand{\fB}{\mathfrak{B}}
\newcommand{\fC}{\mathfrak{C}}

\newcommand{\fF}{\mathfrak{F}}

\newcommand{\fL}{\mathfrak{L}}

\newcommand{\fN}{\mathfrak{N}}

\newcommand{\fP}{\mathfrak{P}}

\newcommand{\fT}{\mathfrak{T}}

\newcommand{\pA}{\mathcal{A}}
\newcommand{\pB}{\mathcal{B}}

\newcommand{\pR}{\mathcal{R}}

\newcommand{\pT}{\mathcal{T}}

\newcommand{\eH}{\mathscr{H}}

\newcommand{\eK}{\mathscr{K}}
\newcommand{\eL}{\mathscr{L}}
\newcommand{\eM}{\mathscr{M}}
\newcommand{\eN}{\mathscr{N}}
\newcommand{\eP}{\mathscr{P}}

\newcommand{\eU}{\mathscr{U}}

\newcommand{\eX}{\mathscr{X}}

\newcommand{\eY}{\mathscr{Y}}

\newcommand{\bC}{\mathbb{C}}
\newcommand{\bD}{\mathbb{D}}

\newcommand{\bN}{\mathbb{N}}
\newcommand{\bT}{\mathbb{T}}

\newcommand{\osum}{\mathbin{\rotatebox[origin=c]{90}{{\small ${\bowtie}\!$}}}}

\DeclareMathOperator{\Ker}{{\rm ker}}

\DeclareMathOperator{\Rk1}{{\rm Rk_1}}

\DeclareMathOperator{\Alg}{Alg}

\DeclareMathOperator{\Lat}{Lat}
\DeclareMathOperator{\Pin}{Pin}

\newtheorem{theorem}{Theorem}
\newtheorem{proposition}[theorem]{Proposition}
\newtheorem{lemma}[theorem]{Lemma}
\newtheorem{corollary}[theorem]{Corollary}

\theoremstyle{definition}

\numberwithin{equation}{section}

\begin{document}

\title[Reflexive subspace lattices in Banach spaces]{Reflexive subspace lattices in Banach spaces}
\author[J. Bra\v{c}i\v{c}]{Janko Bra\v{c}i\v{c}}
\address{Faculty of Natural Sciences and Engineering, University of Ljubljana, A\v{s}ker\v{c}eva c. 12, SI-1000 Ljubljana, Slovenia}
\email{janko.bracic@ntf.uni-lj.si}


\begin{abstract}
In this survey paper we present known results about reflexive subspace lattices. We show that
every nest and every atomic Boolean subspace lattice in a complex Banach space is reflexive, even 
strongly reflexive. Our main tool is Ringrose's Lemma about presence of rank-one operators in a reflexive 
algebra of operators. We consider realizations of small lattices as subspace lattices in Banach spaces.
We exhibit some realizations of the pentagon and prove that some of these realizations are reflexive. 
On the other hand, we show that double triangle subspace lattices in finite-dimensional Banach spaces 
are non-reflexive.
\end{abstract}
\maketitle

\section{Introduction} \label{Sec1}
\setcounter{theorem}{0}

More than fifty years ago, Halmos \cite{Hal} initiated the study of subspace lattices. A few years earlier,
Ringrose \cite{Rin1, Rin2} introduced nest algebras and one of his results says, in the language introduced by
Halmos, that every complete nest in a separable complex Hilbert space is reflexive. Halmos proved that the
same holds for atomic Boolean subspace lattices and he raised several questions related to the reflexivity 
and transitivity of subspace lattices, see \cite{Hal1, Hal, Hal2}. Afterward, many authors studied subspace lattices,
mostly in complex Hilbert spaces, in connection with the theory of non-selfadjoint weakly closed algebras of operators. 
For instance, in \cite{Arv} Arveson developed some difficult machinery to prove that every strongly closed 
commutative subspace lattice is reflexive (commutative means that the projections onto subspaces in the 
lattice commute with each other). This is a generalization of Ringrose's Theorem for nests. Later Davidson
\cite{Dav}, Shulman \cite{Shu}, and Haydon \cite{Hay} gave three different and somewhat simpler proofs of 
Arveson's theorem. 

In this survey paper, we will not follow Arveson but we will rather focus on elementary results 
that can be formulated and proved for general Banach spaces. This also means that we will avoid results obtained 
by Harrison and Longstaff about reflexivity (and strong reflexivity) of certain distributive subspace lattices 
\cite{Har, HL, Lon1}. See \cite{Joh} for a proof that each distributive subspace lattice in a finite-dimensional
Banach space is reflexive.

We aim to present basic results about reflexive subspace lattices in Banach spaces. Our exposition is
as simple as possible and the included proofs are elementary. For the proofs, we do not claim any originality as 
they are 
mostly suitable modifications of the proofs which can be found in the cited literature. It seems that at this level
the main tool is the set of rank-one operators in the algebra induced by a subspace lattice. Therefore, we begin 
our discussion in \Cref{Sec3} by proving a quite general version of Ringrose's Lemma about the presence of 
rank-one operators in reflexive algebras. It was observed by many authors that rank-one and finite rank operators 
play a crucial role in the study of these algebras, see \cite{Erd, HM, Lam1, LauL, Lon3}. \Cref{Sec4} is devoted
to nests. We follow ideas from \cite{Rin1, Rin2} and prove that every complete nest in a Banach space is reflexive.
Moreover, we give direct proof (without using complete distributivity) that nests are strongly reflexive.
Similar results are obtained in \Cref{Sec5} for atomic Boolean algebras. We prove in \Cref{Sec6} that the ordinal 
sum of two reflexive subspace lattices is reflexive; a result that was known to Halmos \cite{Hal}. This enables to build
reflexive subspace lattices from smaller ones. In the last section, we consider the two most important small
lattices, the pentagon and the double triangle. A pentagon subspace lattice can be refelexive or non-reflexive,
while for the double triangle subspace lattices no reflexive example is known. Pentagon and double triangle subspace
lattices have been studied in \cite{DWY, HLR, KLL, LR}.

\section{Preliminaries} \label{Sec2}
\setcounter{theorem}{0}

Let $\eX$ be a complex Banach space of dimension at least $2$. With $\pB(\eX)$ we will denote the Banach
algebra of all operators, i.e., bounded linear transformations, on $\eX$; more generally, $\pB(\eX,\eY)$ is
the Banach space of all operators from a Banach space $\eX$ to Banach space $\eY$.
A subspace of $\eX$ means a closed linear manifold. Subspaces $\{ 0\}$ and $\eX$ are said to be trivial.
If $\eU\subseteq \eX$ is a non-empty set of vectors, then we denote by $[\eU]$ the smallest subspace of $\eX$
that contains $\eU$, that is, $[\eU]$ is the closed linear span of $\eU$. 

The dual space of $\eX$ is denoted by $\eX^*$ and the pairing between $\eX$ and $\eX^*$ is given by 
$\langle x,\xi\rangle=\xi(x)$, for all $x\in \eX$ and $\xi\in \eX^*$.
The annihilator of a subspace $\eM\subseteq \eX$ is the subspace
$\eM^\perp=\{ \xi\in \eX^*;\; \langle x,\xi\rangle=0,\; \text{for all}\; x\in \eM\}$ of $\eX^*$. 
By the Hahn-Banach Theorem, $\eM=\bigcap\limits_{\xi\in \eM^\perp}\Ker(\xi)$,
where $\Ker(\xi)=\{ x\in \eX;\; \langle x,\xi\rangle=0\}$ is the kernel of $\xi$. 
It is obvious that $\eM\subseteq \eN$ if and only if $\eM^\perp\supseteq \eN^\perp$, for
arbitrary subspaces $\eM$ and $\eN$.

Let $\fC(\eX)$ be the family of all subspaces of $\eX$. It is obvious that, for $\eM,\eN\in \fC(\eX)$, the join
$\eM\vee \eN=\overline{\eM+\eN}$ (the bar denotes the closure) and the meet $\eM\wedge\eN=\eM\cap\eN$ 
are in $\fC(\eX)$. Hence, $\fC(\eX)$ is a lattice. A family $\fL\subseteq \fC(\eX)$ is a sublattice
if $\eM\vee\eN$ and $\eM\wedge\eN$ are in $\fL$ whenever $\eM$ and $\eN$ are in $\fL$. A sublattice
$\fL\subseteq \fC(\eX)$ is complete if, for an arbitrary family $\{ \eM_j;\;j\in J\}\subseteq \fL$,
subspaces $\bigvee\limits_{j\in J}\eM_j$ and $\bigwedge\limits_{j\in J}\eM_j$ are in $\fL$.
Here we adopted the convention that $\bigvee \emptyset=\{ 0\}$ and $\bigwedge\emptyset=\eX$.
It is clear that $\fC(\eX)$ is a complete lattice.
A complete sublattice $\fL\subseteq \fC(\eX)$ is called a {\em subspace lattice} if it contains the trivial
subspaces $\{ 0\}$ and $\eX$.

By an abstract lattice we mean an algebraic structure $(L,\vee,\wedge)$, consisting of a non-empty set $L$
and two binary operations that are commutative and associative and satisfy the absorption laws 
$a\vee(a\wedge b)=a$ and $a\wedge(a\vee b)=a$, for all $a, b\in L$. One can define a partial order $\leq$ on $L$ 
by setting $a\leq b$ if $a=a\wedge b$ or, equivalently, $b=a\vee b$. We will always assume that $L$  has a 
greatest element $1$ and a least element $0$. Hence, $0\leq a\leq 1$, for all $a\in L$. If $a, b\in L$ are such that
$a\leq b$, then $[a,b]=\{ c\in L;\; a\leq c\leq b\}$ is an {\em interval} in $L$.

An isomorphism of lattices 
$L_1$ and $L_2$ is a bijective mapping $\phi\colon L_1\to L_2$ such that $\phi(a\vee b)=\phi(a)\vee\phi(b)$,
$\phi(a\wedge b)=\phi(a)\wedge\phi(b)$, and $\phi(1)=1$, $\phi(0)=0$, equivalently, $a\leq b$ (in $L_1$) 
if and only if $\phi(a)\leq \phi(b)$ (in $L_2$). A bijective mapping $\psi\colon L_1\to L_2$ such that 
$\psi(a\vee b)=\psi(a)\wedge\psi(b)$, $\psi(a\wedge b)=\psi(a)\vee\psi(b)$, and $\psi(1)=0$, $\psi(0)=1$
is an anti-isomorphism of lattices; hence, $a\leq b$ (in $L_1$) if and only if $\psi(a)\geq \psi(b)$ (in $L_2$). 
An abstract lattice $L$ is {\em realizable} in a Banach space $\eX$ if there exists a lattice of subspaces 
$\fL\subseteq \fC(\eX)$ and an isomorphism $\phi\colon L\to \fL$. 
In this case, we say that $\fL$ is a realization of $L$. 

Let $\eX$ be a Banach space and let $n\in \bN$. Denote by $\eX^{(n)}$ the direct sum of $n$-copies of $\eX$
equipped with the norm $\| (x_1,\ldots,x_n)\|=\sqrt{\| x_1\|^2+\cdots+\| x_n\|^2}$. Of course, each $\eX^{(n)}$
is a Banach space.
Let $\eX^{(\infty)}$ be the vector space of all sequences $(x_j)_{j=1}^{\infty}\subseteq \eX$ such that
$\sum\limits_{j=1}^{\infty}\| x_j\|^2<\infty$. It is well-known and not hard to see that $\eX^{(\infty)}$
is a Banach space when it is equipped with the norm 
$\| (x_j)_{j=1}^{\infty}\|=\bigl( \sum\limits_{j=1}^{\infty}\| x_j\|^2\bigr)^{1/2}$.
For $\eM\in \fC(\eX)$, let $\eM^{(n)}=\bigl\{ (x_1,\ldots,x_n)\in \eX^{(n)};\; x_j\in \eM,\; j=1,\ldots,n\bigr\}$.
It is obvious that $\eM^{(n)}\in \fC(\eX^{(n)})$. Similarly, one can show that
$\eM^{(\infty)}=\bigl\{ (x_j)_{j=1}^{\infty}\in \eX^{(\infty)};\; x_j\in \eM,\; j\in \bN\bigr\}$
is in $\fC(\eX^{(\infty)})$.

\begin{lemma} \label{lem13}
Let $\fL\subseteq \fC(\eX)$ be a realization of $L$ in a Banach space $\eX$. Then, for every $n\in\bN\cup\{\infty\}$,
$\fL^{(n)}=\{ \eM^{(n)};\; \eM\in \fL\}$ is a realization of $L$ in $\eX^{(n)}$.
\end{lemma}

\begin{proof}
Let $\phi\colon L\to \fL$ be an isomorphism of lattices. Hence, $\fL=\{ \phi(a);\; a\in L\}$
and $a\leq b$ if and only if $\phi(a)\subseteq \phi(b)$. Let $n\in \bN\cup\{\infty\}$. For every $a\in L$,
define $\phi^{(n)}(a)=\phi(a)^{(n)}\in \fC(\eX^{(n)})$. It is obvious that $\phi^{(n)}\colon L\to \fL^{(n)}$
is a bijection such that $a\leq b$ if and only if $\phi^{(n)}(a)\subseteq \phi^{(n)}(b)$.
\end{proof}

Let $L$ be a lattice. Element $a$ is said to be an atom if $0\ne a\ne 1$ and $0\leq b\leq a$, for some $b\in L$,
implies that either $b=0$ or $b=a$. By $L_a$ we will denote the set of all atoms in $L$ (of course, it can be empty).
Lattice $L$ is said to be 
\begin{itemize}
\item[(i)] {\em complemented} if for every $a\in L$ there exists $a'\in L$ such that $a\wedge a'= 0$ and 
$a\vee a'=1$;
\item[(ii)] {\em modular} if, for $a,b,c\in L$, it follows from $a\leq b$ that
$a\vee(c\wedge b)=(a\vee c)\wedge b$;
\item[(iii)] {\em distributive} if equality
$ a\wedge(c\vee b)=(a\wedge c)\vee(a\wedge b)$
and its dual form
$a\vee(c\wedge b)=(a\vee c)\wedge(a\vee b)$
hold for all $a, b, c\in L$;
\item[(iv)] {\em atomic} if the family of all atoms $L_a$ generates $L$ in the sense
that $b=\bigvee\{a\in L_a;\;a\leq b\}$, for every $b\in L$;
\item[(v)] {\em Boolean} if it is complemented and distributive;
\item[(vi)] a {\em nest} if it is totally ordered (i.e., if $a, b \in \fL$, then either $a\leq b$
or $b\leq a$).
\end{itemize}
Every distributive lattice is modular. Indeed, if we assume that $a\leq b$, then the first equality in (iii) 
simplifies to the equality in (ii).

\begin{lemma} \label{lem06}
Every subspace lattice in a finite-dimensional Banach space is modular.
\end{lemma}

\begin{proof}
Let $\eX$ be a finite-dimensional Banach space. Hence, $\eM\vee\eN=\eM+\eN$ and $\eM\wedge\eN=\eM\cap\eN$, 
for arbitrary subspaces $\eM,\eN\in\fC(\eX)$. Let $\fL\subseteq \fC(\eX)$ be a subspace lattice and
let $\eM,\eN,\eK\in \fL$ be such that $\eM\subseteq\eN$. We have to see that 
$\eM+(\eK\cap\eN)=(\eM+\eK)\cap\eN$. Let $u\in \eM+(\eK\cap\eN)$ be arbitrary. Then there exist
$x\in\eM$ and $y\in \eK\cap\eN$ such that $u=x+y$. Since $x\in\eM\subseteq\eN$ we see that
$x\in (\eM+\eK)\cap\eN$. We also have $y\in (\eM+\eK)\cap\eN$ and therefore $u\in(\eM+\eK)\cap\eN$. 
This proves the inclusion $\eM+(\eK\cap\eN)\subseteq (\eM+\eK)\cap\eN$. To prove the opposite inclusion,
assume that $u\in (\eM+\eK)\cap\eN$. Then $u\in\eN$ and there exist $x\in\eM\subseteq \eN$ and $y\in \eK$
such that $u=x+y$ which gives $y=u-x\in\eN$, more precisely, $y\in\eK\cap\eN$. Hence, 
$u=x+y\in \eM+(\eK\cap\eN)$. We conclude that $\fL$ is modular.
\end{proof}

Let $L$ be a complemented lattice. It is not hard to see that the complement of every $a\in L$ is unique if
$L$ is distributive. In particular, if $L$ is a Boolean lattice then every element has a unique complement.
Modularity is not enough for the uniqueness of a complement. Indeed, by \Cref{lem06}, every subspace
lattice in a finite-dimensional Banach space $\eX$ is modular, but it is clear that a non-trivial subspace 
$\eM\subseteq \eX$ does not have a unique complement. Here we want to stress the difference between
the notion of a complement in the context of subspace lattices and the notion of a complement of
a subspace in the theory of Banach spaces. In the theory of Banach spaces, a subspace $\eM\subseteq \eX$
is complemented if it is algebraically complemented, that is, there exists a subspace $\eN\subseteq \eX$ 
such that  $\eM\cap\eN=\{ 0\}$ and $\eM+\eN=\eX$. If $\eM,\eN\in \fC(\eX)$ are such that
$\eM\cap\eN=\{ 0\}$ and $\eM+\eN$ is dense in $\eX$, then it is said that $\eN$ is a {\em quasi-complement}
of $\eM$; we say that $(\eM,\eN)$ is a {\em quasi-complemented pair} of subspaces. 
In the context of subspace lattices,
$\eM$ has a complement (or, is complemented) in a subspace lattice $\fL$ if there exists $\eN\in \fL$ such that 
$\eM\wedge\eN=\{ 0\}$ and 
$\eM\vee\eN=\eX$, that is, $\eM\cap\eN=\{ 0\}$ and $\overline{\eM+\eN}=\eX$ (i.e., $(\eM,\eN)$ is a 
quasi-complement pair of subspaces.
Of course, if $\eM$ is finite-dimensional or has a finite codimension, then both notions coincide.

A subspace $\eM\in \fC(\eX)$ is invariant for an operator $T\in \pB(\eX)$ if $T\eM\subseteq \eM$.
More generally, a subspace $\eM$ is invariant for a non-empty set of operators $\pT$, if $\eM$
is invariant for every operator in $\pT$. Let $\Lat(\pT)$ denote the family of all subspaces which are
invariant for $\pT$. It is not hard to see that $\Lat(\pT)$ is a subspace lattice.

For a non-empty family $\fF\subseteq \fC(\eX)$, let 
$\Alg(\fF)=\{ A\in \pB(\eX);\; A\eM\subseteq \eM,\; \text{for all}\; \eM\in \fF\}$. This is a strongly closed
unital subalgebra of $\pB(\eX)$. The following is well-known and easy to see: $\pT\subseteq \Alg\Lat(\pT)$,
$\Lat(\pT)=\Lat\Alg\Lat(\pT)$, $\fF\subseteq \Lat\Alg(\fF)$, and $\Alg(\fF)=\Alg\Lat\Alg(\fF)$.

For a subspace lattice $\fL\subseteq \fC(\eX)$ we have $\fL\subseteq \Lat\Alg(\fL)\subseteq \fC(\eX)$.
Of course, lattices for which one of these inclusions is equality are the most interesting. In this paper
we are concerned with {\em reflexive} subspace lattices, that is, those satisfying $\Lat\Alg(\fL)=\fL$.
This definition is due to Halmos \cite{Hal}. Lattices $\fL$ satisfying the other extreme condition, i.e., 
$\Lat\Alg(\fL)= \fC(\eX)$ are said to be {\em transitive}. 
A unital strongly closed subalgebra $\pA\subseteq \pB(\eX)$ is reflexive if $\Alg\Lat(\pA)=\pA$. Hence,
a subspace lattice $\fL$ is reflexive if and only if $\fL=\Lat(\pT)$, for a subset $\pT\subseteq \pB(\eX)$,
and a strongly closed algebra $\pA$ is reflexive if and only if $\pA=\Alg(\fF)$, for a family $\fF\subseteq \fC(\eX)$.

We have already mentioned that we are interested in the reflexive subspace lattices. At this point, 
we give a simple characterization of these lattices. We need the notion of a cyclic subspace. Let 
$\pA\subseteq \pB(\eX)$ be an algebra. The cyclic subspace of $\pA$ generated by $x\in \eX$ is $[\pA x]$, that is,
the closed linear span of $\{ Ax;\; A\in \pA\}$. It is obvious that $[\pA x]\in \Lat(\pA)$, for all $x\in \eX$.

\begin{proposition} \label{prop04}
For a subspace lattice $\fL\subseteq \fC(\eX)$, the following are equivalent:
\begin{itemize}
\item[(a)] $\fL$ is reflexive;
\item[(b)] $\bigl[ \Alg(\fL)x\bigr]\in \fL$, for all $x\in \eX$;
\item[(c)] there exists a subset $\pR$ of $\Alg(\fL)$ such that $\fL=\Lat(\pR)$.
\end{itemize}
\end{proposition}

\begin{proof}
(a)$\iff$(b). If $\fL$ is reflexive, then $\bigl[ \Alg(\fL)x\bigr]\in \Lat\Alg(\fL)=\fL$, for all $x\in \eX$. 
To prove the opposite implication, assume that $\eM\in \Lat\Alg(\fL)$.
Since $\eM=\bigvee\limits_{x\in \eM}\bigl[ \Alg(\fL)x\bigr]$ and $\bigl[ \Alg(\fL)x\bigr]\in \fL$, for all $x\in \eM$,
the completeness of $\fL$ gives $\eM\in \fL$.

(a)$\iff$(c). If $\fL$ is reflexive, then we can take $\pR=\Alg(\fL)$. On the other hand, it is clear that
$\fL=\Lat(\pR)$ implies $\Lat\Alg(\fL)=\Lat(\pR)=\fL$.
\end{proof}

As we shall see later (see \Cref{theo06} and \Cref{theo08}), the case when $\pR$ in the assertion (c) is the set 
of all rank-one operators in $\Alg(\fL)$ is of particular interest.

\section{Rank-one operators in $\Alg(\fL)$} \label{Sec3}
\setcounter{theorem}{0}

In this section, we will present Ringrose's result \cite[Lemma 3.3]{Rin1} about the presence of rank-one operators 
in reflexive algebras. Ringrose proved his lemma for nest algebras on Hilbert spaces.
Later Longstaff \cite[Lemma 3.1]{Lon1} proved it for arbitrary reflexive algebras on Hilbert spaces.
Lambrou \cite[Lemma 2.1]{Lam} observed that the lemma holds for every reflexive algebra on
complex Banach spaces. Before we state and prove the lemma we need to introduce new notions.

If $\fL\subseteq \fC(\eX)$ is a subspace lattice, then, for an arbitrary subspace $\eM\in \fC(\eX)$, let 
\begin{equation} \label{eq01}
\eM_{-}=\bigvee\{ \eK\in \fL;\; \eM\not\subseteq \eK\}\qquad \text{and}\qquad
\eM_{+}=\bigwedge\{ \eK\in \fL;\; \eK\not\subseteq \eM\}.
\end{equation}
Since $\fL$ is a complete lattice we see that $\eM_{-}$ and $\eM_{+}$ are subspaces in $\fL$. 
It is easily seen that $\{ 0\}_{-}=\{0\}$, $\{ 0\}_{+}=\bigwedge\limits_{\substack{\eK\in\fL\\ \eK\ne \{0\}}}\eK$,
$\eX_{-}=\bigvee_{\substack{\eK\in\fL\\ \eK\ne \eX}}\eK$, and $\eX_{+}=\eX$.

For $0\ne e\in\eX$ and $0\ne \eta\in \eX^*$, let $e\otimes\eta\colon \eX\to\eX$ be the rank-one operator 
that is given by $(e\otimes\eta)x=\langle x,\eta\rangle e$, for all $x\in \eX$. 
For an algebra $\pA\subseteq \pB(\eX)$, let $\Rk1(\pA)$ be the set of all rank-one operators in $\pA$. 
Note that the linear span of $\Rk1(\pA)$ is a  subalgebra of $\pA$ if $\Rk1(\pA)$ is not empty.

\begin{lemma}[Ringrose] \label{lem01}
Let $\fL\subseteq \fC(\eX)$ be a subspace lattice. A rank-one operator $e\otimes\eta$ belongs to the algebra
$\Alg(\fL)$ if and only if there is a subspace $\eM\in \fL$ such that $e\in \eM$ and $\eta\in \eM_{-}^{\perp}$.
\end{lemma}

\begin{proof}
Let $\eM\in \fL$ and let $e\in \eM$, $\eta\in \eM_{-}^{\perp}$ be arbitrary. If $\eN\in \fL$, then either 
$\eM\subseteq \eN$ or $\eM\not\subseteq \eN$. In the former case, $(e\otimes \eta)x=\langle x,\eta\rangle e\in \eN$,
for all $x\in \eN$, and in the latter case, $(e\otimes \eta)x=\langle x,\eta\rangle e=0$,
for all $x\in \eN$, since $\eN\subseteq \eM_{-}$ and therefore $\eN^\perp\supseteq \eM_{-}^{\perp}$. 
In any case, $\eN$ is invariant for $e\otimes\eta$ and therefore $e\otimes\eta\in \Alg(\fL)$.

Conversely, let $0\ne e\otimes \eta\in \Alg(\fL)$. Let $\eM=\bigcap\{ \eN\in \fL;\; e\in \eN\}$. It is clear that
$\eM\in \fL$ and $e\in \eM$. Suppose that $\eN\in \fL$ is such that $\eM\not\subseteq \eN$. It follows, by
the definition of $\eM$, that $e\not\in \eN$. Since $\eN$ is invariant for $e\otimes\eta$ it follows from
$(e\otimes \eta)x=\langle x,\eta\rangle e$ that $\langle x,\eta\rangle=0$, for all $x\in \eN$. Hence,
$\eN\subseteq \Ker(\eta)$, for all subspaces $\eN\in \fL$ such that $\eM\not\subseteq \eN$,
and therefore $\eM_{-}\subseteq \Ker(\eta)$, that is, $\eta\in \eM_{-}^{\perp}$.
\end{proof}

Since $\Alg(\fL)=\Alg\Lat\Alg(\fL)$ it follows from \Cref{lem01} that a rank-one operator $e\otimes \eta$
belongs to $\Alg(\fL)$ if and only if there exists a subspace $\eM\in \Lat\Alg(\fL)$ such that $e\in \eM$
and $\eta\in \eM_{-}^{\perp}$.

\begin{corollary} \label{cor03}
Let $\fL\subseteq \fC(\eX)$ be a subspace lattice and let $e\in \eX$ be arbitrary. A rank-one operator $e\otimes \eta$
belongs to $\Alg(\fL)$ if and only if $\eta \in \bigl[ \Alg(\fL)e\bigr]_{-}^{\perp}$.
\end{corollary}

\begin{proof}
Since $\bigl[ \Alg(\fL)e\bigr]\in \Lat\Alg(\fL)$ and $e\in \bigl[ \Alg(\fL)e\bigr]$ it follows, by \Cref{lem01},
that $e\otimes\eta\in \Alg(\fL)$ whenever $\eta\in \bigl[ \Alg(\fL)e\bigr]_{-}^{\perp}$.

Conversely, assume that $e\otimes\eta\in \Alg(\fL)$. If $\eK\in \Lat\Alg(\fL)$ is such that 
$\bigl[ \Alg(\fL)e\bigr]\not\subseteq \eK$, then there exists $T\in \Alg(\fL)$ for which $Te\not\in \eK$.
Since $\eK$ is invariant for $T$ it follows that $e\not\in \eK$. Suppose that there
existed $f\in \eK$ such that $\langle f,\eta\rangle \ne 0$. It would follow that
$\langle f,\eta\rangle e=(e\otimes \eta)f\in \eK$ because $e\otimes \eta\in \Alg(\fL)$. But this is not the case
since $e\not\in\eK$.
Hence, $\eK\subseteq \Ker(\eta)$ and consequently $\bigl[ \Alg(\fL)e\bigr]_{-}\subseteq \Ker(\eta)$.
\end{proof}

Since $\Rk1(\Alg(\fL))\subseteq \Alg(\fL)$ we have $\Lat\Alg(\fL)\subseteq \Lat\bigl(\Rk1(\Alg(\fL))\bigr)$. 
The following theorem which is based on \cite[Proposition 3.2]{Lon1} characterizes subspaces in 
$\Lat\bigl(\Rk1(\Alg(\fL))\bigr)$.

For $\eN\in \fC(\eX)$, let
\begin{equation} \label{eq04}
\eN_{*}=\bigwedge\{ \eM_{-};\; \eM\in \fL\;\;\text{and}\;\;\eM\not\subseteq \eN\}.
\end{equation}
It is clear that $\eN_{*}\in\fL$ and it is easily seen that 
$\{ 0\}_{*}=\bigwedge\limits_{\substack{\eM\in\fL\\ \eM\ne \{0\}}}\eM_{-}$ and $\eX_{*}=\eX$. 

\begin{theorem} \label{theo07}
Let $\fL\subseteq \fC(\eX)$ be a subspace lattice. For $\eK\in \fC(\eX)$, the following assertions are equivalent:
\begin{itemize}
\item[(a)] $\eK\in \Lat\bigl(\Rk1(\Alg(\fL))\bigr)$;
\item[(b)] for every $\eM\in \fL$, either $\eM\subseteq \eK$ or $\eK\subseteq \eM_{-}$;
\item[(c)] there exists a subspace $\eN\in \fL$ such that $\eN\subseteq \eK\subseteq \eN_{*}$.
\end{itemize}
\end{theorem}

\begin{proof}
(a)$\Rightarrow$(b). Let $\eK\in \Lat\bigl(\Rk1(\Alg(\fL))\bigr)$ and let $\eM\in \fL$ be arbitrary. If $\eM_{-}=\eX$,
then there is nothing to prove, that is, $\eK\subseteq \eM_{-}$. Assume therefore that $\eM_{-}\ne \eX$.
There are two possibilities: either $\eK\subseteq \ker(\eta)$, for all $\eta\in \eM_{-}^{\perp}$, or there
exists $\eta\in \eM_{-}^{\perp}$ such that $\eK\not\subseteq \Ker(\eta)$. The former case gives
$\eM_{-}^{\perp}\subseteq\eK^{\perp}$ and therefore $\eK\subseteq \eM_{-}$. In the latter case, there
exists $y\in \eK$ such that $\langle y,\eta\rangle=1$. Since $e\otimes \eta\in \Rk1(\Alg(\fL))$, for all $e\in \eM$,
by \Cref{lem01}, we have $(e\otimes\eta)y=e\in \eK$, that is, $\eM\subseteq \eK$.

(b)$\Rightarrow$(c).
Let $\eK\in \fC(\eX)$ be arbitrary. Denote $\eN=\bigvee\{\eM\in \fL;\; \eM\subseteq \eK\}$. It is obvious
that $\eN\in \fL$ and $\eN\subseteq \eK$. Hence, if $\eM\in \fL$ is such that $\eM\not\subseteq \eN$,
then $\eM\not\subseteq \eK$. Thus, if $\eM\in \fL$ and $\eM\not\subseteq \eN$, then $\eK\subseteq \eM_{-}$,
by (b). It follows that 
$\eK\subseteq \bigwedge\{ \eM_{-};\; \eM\in \fL\;\;\text{and}\;\;\eM\not\subseteq \eN\}=\eN_{*}$.

(c)$\Rightarrow$(a).
Let $\eK\in \fC(\eX)$ be arbitrary and suppose that there exists $\eN\in \fL$ such that
$\eN\subseteq \eK\subseteq \eN_{*}$. Let $e\otimes \eta$ be an arbitrary rank-one operator in $\Alg(\fL)$.
By \Cref{lem01}, there exists $\eM\in \fL$ such that $e\in \eM$ and $\eta\in \eM_{-}^{\perp}$.
If $\eM\subseteq \eN$, then $(e\otimes\eta)x=\langle x,\eta\rangle e\in \eM\subseteq \eN\subseteq \eK$,
for all $x\in \eK$. On the other hand, if $\eM\not\subseteq \eN$, then $\eK\subseteq \eM_{-}$, by the
assumption. Hence, $(e\otimes\eta)x=0\in \eK$, for all $x\in \eK$. We have seen that $\eK$ is
invariant for every rank-one operator in $\Alg(\fL)$ and therefore we may conclude that 
$\eK\in \Lat\bigl(\Rk1(\Alg(\fL))\bigr)$.
\end{proof}

\begin{corollary} \label{cor04}
Let $\fL\subseteq\fC(\eX)$ be a subspace lattice. If $\eN_{*}=\eN$, for all $\eN\in \fL$, then $\fL$ is reflexive;
moreover, $\fL=\Lat\bigl(\Rk1(\Alg(\fL))\bigr)$.
\end{corollary}

\begin{proof}
Since $\Rk1(\Alg(\fL))\subseteq \Alg(\fL)$ we have 
$\fL\subseteq \Lat\Alg(\fL)\subseteq \Lat\bigl(\Rk1(\Alg(\fL))\bigr)$. On the other hand, if 
$\eK\in \Lat\bigl(\Rk1(\Alg(\fL))\bigr)$, then there exists $\eN\in \fL$ such that 
$\eN\subseteq \eK\subseteq \eN_{*}$, by \Cref{theo07}. However, by the assumption, 
$\eN_{*}=\eN$ and therefore $\eK=\eN\in \fL$.
\end{proof}

A subspace lattice $\fL\subseteq\fC(\eX)$ such that $\eN_{*}=\eN$, for all $\eN\in \fL$, is said to be
{\em strongly reflexive} (see, \cite[Definition 4.2]{Lon1}). Hence, by \Cref{cor04}, strongly reflexive
subspace lattices are reflexive. Longstaff has proved that a subspace lattice $\fL\in \fC(\eH)$, where $\eH$ is 
a complex Hilbert space, is strongly reflexive if and only if it is {\em completely distributive}. We will not
use complete distributivity, so the reader is referred to \cite{Lon1} for its precise definition. Nests and 
atomic Boolean subspace lattices are completely distributive and therefore (in the case of Hilbert space)
are strongly reflexive, according to Longstaff's result. We will give direct proof of the 
strong reflexivity of nests and atomic Boolean subspace lattices (on an arbitrary Banach space) without
using the fact that these subspace lattices are completely distributive, see \Cref{theo06} and \Cref{theo08}.

\section{Reflexivity of nests} \label{Sec4}
\setcounter{theorem}{0}

In what follows, nest means a totally ordered subspace lattice. In this section, we will show that nests are 
strongly reflexive. First, we will follow the ideas from Ringrose's papers \cite{Rin1} and \cite{Rin2} and 
show that every nest is reflexive. At the end of the section, we will show that every nest is strongly reflexive.

A subspace $\eP\in \fL$ is a {\em pinch point} of $\fL$ if, for every $\eM\in \fL$, either $\eM\subseteq \eP$ or 
$\eP\subseteq \eM$ (Halmos \cite{Hal} call pinch points comparable subspaces). 
Let $\Pin(\fL)$ be the family of
all pinch points of $\fL$. Of course, a subspace lattice $\fL$ is a nest if and only if $\Pin(\fL)=\fL$.

\begin{lemma} \label{lem07}
If $\fL$ is a subspace lattice, then $\Pin(\fL)$ is a subspace lattice and it is a nest.
\end{lemma}

\begin{proof}
It is clear that $\Pin(\fL)$ is a totally ordered sublattice of $\fL$ containing $\{ 0\}$ and $\eX$. To see that it is
complete, let $\{ \eP_j;\;j\in J\}\subseteq \Pin(\fL)$ be an arbitrary non-empty family. If $\eM\in \fL$ is such that
$\eM\subseteq \eP_j$, for all $j\in J$, then $\eM\subseteq \bigwedge_{j\in J}\eP_j\subseteq
\bigvee_{j\in J}\eP_j$. Similarly, if $\eM\supseteq \eP_j$, for all $j\in J$, then 
$\eM\supseteq \bigvee_{j\in J}\eP_j\supseteq \bigwedge_{j\in J}\eP_j$. Assume that there exist indices $j_1, j_2\in J$ such that $\eM\subseteq \eP_{j_1}$ and $\eM\supseteq \eP_{j_2}$. Then $\eM\subseteq\bigvee_{j\in J}\eP_j$
and $\eM\supseteq\bigwedge_{j\in J}\eP_j$. Thus, in any case, $\bigvee_{j\in J}\eP_j$ and 
$\bigwedge_{j\in J}\eP_j$ are pinch points of $\fL$.
\end{proof}

\begin{lemma} \label{lem05}
If $\eP\in \fL$ is a pinch point, then
\begin{equation} \label{eq02}
\eP_{-}=\bigvee\{ \eK\in \fL;\; \eK\subsetneq \eP\}\subseteq \eP\quad \text{and}\quad
\eP_{+}=\bigwedge\{ \eK\in \fL;\; \eP\subsetneq \eK\}\supseteq \eP.
\end{equation}
Moreover, if $\eM\in \fL$ is such that $\eP_{-}\subseteq \eM\subseteq \eP_{+}$, then $\eM=\eP_{-}$,
or $\eM=\eP$, or $\eM=\eP_{+}$.
\end{lemma}

\begin{proof}
Since $\eP\in \fL$ is a pinch point we have $ \eK\not\subseteq \eP \iff \eP\subsetneq \eK$ and
$\eP\not\subseteq \eK \iff \eK\subsetneq \eP$, for all $\eK\in \fL$. Thus, $\eP_{-}$
is the span of all subspaces from $\fL$ that are properly contained in $\eP$ and, similarly,
$\eP_{+}$ is the intersection of all subspaces from $\fL$ that contain $\eP$ as a proper subspace.
It is clear that this implies $\eP_{-}\subseteq \eP\subseteq \eP_{+}$. Let $\eM\in \fL$ be
such that $\eP_{-}\subseteq \eM\subseteq \eP_{+}$. Suppose that $\eM\ne \eP$. Since $\eP$ is a
pinch point we have either $\eM\subsetneq \eP$ or $\eP\subsetneq \eM$. In the former case, 
$\eM\subseteq \eP_{-}$ and therefore $\eM=\eP_{-}$, since we assumed that 
$\eP_{-}\subseteq \eM\subseteq \eP_{+}$. Similarly, $\eP\subsetneq \eM$ gives $\eM=\eP_{+}$.
\end{proof}

The following lemma is an extension of \cite[Lemma]{Hal} to Banach spaces.

\begin{lemma} \label{lem08}
Let $\fL\subseteq \fC(\eX)$ be a subspace lattice. Then $\Pin(\fL)\subseteq \Pin\bigl(\Lat\Alg(\fL)\bigr)$.
\end{lemma}

\begin{proof}
Let $\eP$ be a pinch point of $\fL$. If $\eK\in \Lat\Alg(\fL)$ is such that $\eK\not\subseteq \eP$,
then there exists $f\in \eK$ such that $f\not\in \eP$ and therefore, by \eqref{eq02}, $f\not\in \eP_{-}$.
Hence, there exists $\eta\in \eP_{-}^{\perp}$ such that $\langle f,\eta\rangle=1$. By \Cref{lem01},
$e\otimes \eta\in \Alg(\fL)$, for all $e\in \eP$. Since $\eK\in \Lat\Alg(\fL)$ it is invariant for $e\otimes \eta$.
Thus, $(e\otimes \eta)f=e\in \eK$, for all $e\in \eP$, that is, $\eP\subseteq \eK$.
\end{proof}

Now we are prepared to generalize Ringrose's result (see \cite[Theorem 3.4]{Rin1}) to general Banach spaces.

\begin{theorem} \label{theo01}
If a subspace lattice $\fN\subseteq \fC(\eX)$ is a nest, then it is reflexive.
\end{theorem}

\begin{proof}
Towards a contradiction, suppose that $\fN$ is not reflexive. Then there exists $\eM\in \Lat\Alg(\fN)$ 
such that $\eM\not\in \fN$. Since every $\eK\in \fN$ is a pinch point for $\Lat\Alg(\fN)$, by \Cref{lem08},
we have either $\eK\subseteq \eM$ or $\eM\subseteq \eK$, for all $\eK\in \fN$. Hence,
$$ \eM_{-}=\bigvee\{ \eK\in \fN;\; \eK\subseteq \eM\}\qquad \text{and}\qquad
\eM_{+}=\bigwedge\{ \eK\in \fN;\; \eM\subseteq \eK\}. $$
It is clear that $\eM_{-}\subseteq \eM\subseteq \eM_{+}$.
Since $\fN$ is a complete lattice we see that $\eM_{-}, \eM_{+}\in \fN$. Hence, 
$\eM_{-}\subsetneq \eM\subsetneq \eM_{+}$ because of $\eM\not\in \fN$.
In particular, $\eM_{-}\ne \eM_{+}$. We distinguish two cases.

{\em Case 1:} There exists $\eL\in \fN$ such that $\eM_{-}\subsetneq \eL\subsetneq \eM_{+}$.\\
If $\eL$ were a subset of $\eM$, then by the definition of $\eM_{-}$ we would have
$\eL\subseteq \eM_{-}$ which is not the case as we have assumed that
$\eM_{-}\subsetneq \eL$. Thus,  $\eL\not\subseteq\eM$ and therefore $\eM^\perp\not\subseteq \eL^\perp$.
Similarly, $\eM\subseteq \eL$ would give $\eM_{+}\subseteq \eL$ which is not true. Hence, 
$\eM\not\subseteq\eL$ and consequently $\eL^\perp\not\subseteq \eM^\perp$. It follows that there
exist a non-zero vector $e\in \eL\setminus \eM$ and a non-zero functional 
$\eta\in \eL^\perp\setminus\eM^\perp$. Since, by \Cref{lem05}, $\eL_{-}\subseteq \eL$ we have
$e\in \eL$ and $\eta\in \eL_{-}^{\perp}$. By \Cref{lem01}, $e\otimes \eta$ is a rank-one operator in $\Alg(\fN)$.
Since $\eta\not\in \eM^\perp$ there exists $x\in \eM$ such that $\langle x,\eta\rangle\ne 0$. Hence,
$(e\otimes \eta)x=\langle x,\eta\rangle e\not\in \eM$, that is, $\eM\not\in \Lat\Alg(\fN)$, which is a contradiction.

{\em Case 2:} There is no $\eL\in \fN$ such that $\eM_{-}\subsetneq \eL\subsetneq \eM_{+}$.\\
It follows that $\eM_{-}$ is the largest space in $\fN$ which is properly contained in $\eM_{+}$.
Hence $\eM_{-}=(\eM_{+})_{-}$. Since $\eM_{-}\subsetneq \eM$ there exists 
$\eta\in (\eM_{+})_{-}^{\perp}$ and $u\in \eM$ such that $\langle u,\eta\rangle=1$. On the other hand, 
because of $\eM\subsetneq \eM_{+}$ there exists $e\in \eM_{+}\setminus\eM$. By \Cref{lem01},
$e\otimes\eta\in \Alg(\fN)$. Since $u\in \eM$ and $(e\otimes\eta)u=e\not\in \eM$
we have a contradiction $\eM\not\in \Lat\Alg(\fN)$, again.
\end{proof}

As we have already mentioned nests are strongly reflexive. We can prove this using lemmas that we prepared.
Recall that $\Rk1(\Alg(\fN))$ is the subset of rank-one operators in $\Alg(\fN)$.

\begin{theorem} \label{theo06}
If a subspace lattice $\fN\subseteq \fC(\eX)$ is a nest, then it is strongly reflexive;
in particular, $\fN=\Lat\bigl(\Rk1(\Alg(\fN))\bigr)$.
\end{theorem}

\begin{proof}
Let $\eN\in \fN$ be arbitrary. Since $\eX_{*}=\eX$ we may assume that $\eN\ne \eX$.
By the definition of $\eN_{*}$, we have
$$\eN_{*}=\bigwedge\{ \eM_{-};\; \eM\in \fN\;\;\text{and}\;\;\eM\not\subseteq \eN\}
=\bigwedge\{ \eM_{-};\; \eM\in \fN\;\;\text{and}\;\;\eN\subsetneq \eM\}.$$
Let $\eK,\eL\in\fN$ be such that $\eK\subseteq\eL$. Then $\eK_{-}=\bigvee\{\eM\in\fN;\; \eM\subsetneq \eK\}
\subseteq\bigvee\{\eM\in\fN;\; \eM\subsetneq \eL\}=\eL_{-}$. 

By \Cref{lem05}, $\eN\subseteq \eN_{+}$. We distinguish two cases: (i) $\eN\subsetneq \eN_{+}$ and
(ii) $\eN=\eN_{+}$. In the case (i), $(\eN_{+})_{-}=\eN$ and $\eN_{+}$ is the smallest subspace in $\fN$ 
which properly contains $\eN$. Hence, 
$\eN_{*}=\bigwedge\{ \eM_{-};\; \eM\in \fN\;\;\text{and}\;\;\eN\subsetneq \eM\}=(\eN_{+})_{-}=\eN$.
Consider now the case (ii). If $\eM\in \fN$ is such that $\eN\subsetneq \eM$, then $\eN\subseteq \eM_{-}$,
by \Cref{lem05}. Hence, $\eN\subseteq \bigwedge\{ \eM_{-};\; \eM\in \fN\;\;\text{and}\;\;\eN\subsetneq \eM\}=
\eN_{*}$. Hence, $\eN_{+}=\eN\subseteq \eN_{*}\subseteq \eN_{+}$, that is, $\eN_{*}=\eN$.
By \Cref{cor04}, strong reflexivity gives $\fN=\Lat\bigl(\Rk1(\Alg(\fN))\bigr)$.
\end{proof}

\section{Reflexivity of Boolean subspace lattices} \label{Sec5}
\setcounter{theorem}{0}

Now we will consider Boolean subspace lattices. Recall that a Boolean subspace lattice is distributive
and every subspace $\eM$ in it has a unique complement $\eM'$. Halmos \cite{Hal} proved that every atomic 
Boolean subspace lattice in a Hilbert space is reflexive. We will extend this to Banach spaces. Moreover,
we will show that every atomic Boolean subspace lattice is strongly reflexive.

Let $\fB\subseteq \fC(\eX)$ be an atomic Boolean subspace lattice. Denote by $\fB_a$ the family of all
atoms in $\fB$. Since $\fB$ is atomic $\fB_a$ cannot be empty. Let $\eM\in \fB$ be non-trivial, that is, 
$\{0\}\ne \eM\ne \eX$. Then, by the definition of an atomic Boolean lattice, 
$\eM=\bigvee\{ \eL\in \fB_a;\; \eL\subseteq \eM\}$. 

\begin{lemma} \label{lem03}
Let $\fB$ be an atomic Boolean subspace lattice. The complement of a non-trivial $\eM\in \fB$ is 
$\eM'=\bigvee\{ \eK\in \fB_a;\; \eK\wedge\eM=\{ 0\}\}$.
\end{lemma}

\begin{proof}
Let $\eK\in \fB_a$. If $\eK\not\subseteq \eM$, then it follows from $\eK\wedge\eM\subseteq \eK$
that $\eK\wedge\eM=\{0\}$. On the other hand, $\eK\wedge\eM=\{0\}$ implies $\eK\not\subseteq \eM$.
Hence, for $\eK\in \fB_a$, relations $\eK\not\subseteq \eM$ and $\eK\wedge\eM=\{0\}$ are equivalent.

Let $\eK\in \fB_a$ be such that $\eK\not\subseteq \eM$, i.e., $\eK\wedge\eM=\{0\}$. 
Again, $\eK\wedge\eM'\subseteq \eK$ implies that either $\eK\wedge\eM'=\{0\}$ or $\eK\subseteq \eM'$. 
If the former were true, then we would have 
$\eK=\eK\wedge(\eM\vee\eM')=(\eK\wedge\eM)\vee(\eK\wedge\eM')=\{0\}$ since $\fB$ is distributive.
Since an atom cannot be trivial we may conclude that $\eK\subseteq \eM'$. It follows that
$\bigvee\{ \eK\in \fB_a;\; \eK\wedge\eM=\{ 0\}\}\subseteq \eM'$. The opposite inclusion is obvious: 
if $\eK\in \fB_a$ is such that $\eK\subseteq \eM'$, then $\eK\wedge\eM=\{ 0\}$.
\end{proof}

\begin{lemma} \label{lem04}
Let $\fB\subseteq \fC(\eX)$ be a Boolean subspace lattice. If $\eK\in\fB_a$, then 
$\eK_{-}=\eK'=\bigvee\{ \eL\in \fB_a;\; \eL\ne \eK\}$. If $\{ 0\}\ne \eM\in \fB\setminus\fB_a$, then
$\eM_{-}=\eX$.
\end{lemma}

\begin{proof}
Since, for $\eL\in \fB_a$, relations $\eL\wedge\eK=\{0\}$ and $\eL\ne\eK$ are equivalent (see the proof 
of \Cref{lem03}) we see that $\eK'=\bigvee\{ \eL\in \fB_a;\; \eL\ne \eK\}$, by \Cref{lem03}.

By the definition, $\eK_{-}=\bigvee\{ \eM\in \fB;\; \eK\not\subseteq \eM\}$. It is clear that
$\bigvee\{ \eM\in \fB;\; \eK\not\subseteq \eM\}=
\bigvee\{ \eM\in \fB;\; \{0\}\ne\eM\ne\eX\;\;\text{and}\;\; \eK\not\subseteq \eM\}$. It follows,
by the first paragraph in the proof of \Cref{lem03}, that $\eK_{-}=\bigvee\{ \eM\in \fB;\; \eK\wedge\eM=\{0\}\}$.
If $\eK\wedge\eM=\{0\}$, then $\eK'=\eK'\vee(\eK\wedge\eM)=(\eK'\vee\eK)\wedge(\eK'\vee\eM)=\eK'\vee\eM$
and therefore $\eM\subseteq\eK'$. Hence, $\eK_{-}=\bigvee\{ \eM\in \fB;\; \eK\wedge\eM=\{0\}\}$.
Since $\fB$ is atomic we have $\eM=\bigvee\{ \eL\in \fB_a;\; \eL\subseteq \eM\}$, for every $\eM\in \fB$.
It follows that
$\bigvee\{ \eM\in \fB;\; \eK\wedge\eM=\{0\}\}=\bigvee\{ \eL\in \fB_a;\; \eK\wedge\eL=\{0\}\}=\eK'$
(the last equality holds by \Cref{lem03}).

Let $\{ 0\}\ne \eM\in \fB\setminus\fB_a$. Since $\eM$ is not an atom we have $\eM\not\subseteq \eK$,
for all $\eK\in \fB_a$. It follows that $\eM_{-}=\bigvee\{\eL\in \fB;\;\eM\not\subseteq \eL\}\supseteq
\bigvee\{ \eK\in \fB_a;\; \eM\not\subseteq \eK\}=\bigvee\limits_{\eK\in \fB_a}\eK=\eX$.
\end{proof}

\begin{lemma} \label{lem02}
Let $\fB\subseteq \fC(\eX)$ be a Boolean subspace lattice. If $\eK$ is an atom in $\fB$ and $\eN\in \Lat\Alg(\fB)$,
then either $\eK\subseteq \eN$ or $\eN\subseteq \eK'$ which is equivalent to $\eK\subseteq \eN$ or 
$\eK\subseteq \eN'$.
\end{lemma}

\begin{proof}
Let $\eK\in \fB$ be an atom and let $\eN\in \Lat\Alg(\fB)$ be arbitrary. Suppose that $\eK\not\subseteq \eN$.
Then there exists $e\in \eK$ such that $e\not\in \eN$. Since $\{ 0\}\ne \eK\ne \eX$, by the definition od atoms, 
we have $\{ 0\}\ne \eK'\ne \eX$, as well, and, by \Cref{lem04}, $\eK_{-}=\eK'$. Let $0\ne \eta\in (\eK')^\perp$ 
be arbitrary. Then
$e\otimes \eta\in \Alg(\fB)$, by \Cref{lem01}. If there existed $f\in \eN$ such that $\langle f,\eta\rangle\ne 0$,
then we would have $\langle f,\eta\rangle e=(e\otimes \eta)f\in \eN$, which is not the case. Hence, $\eN$ is
annihilated by every $\eta\in (\eK')^\perp$, that is, $(\eK')^\perp\subseteq \eN^\perp$ which gives
$\eN\subseteq \eK'$.
\end{proof}

\begin{corollary} \label{cor01}
Let $\fB$ be an atomic Boolean subspace lattice. If $\eM\in \fB$, then 
$$\eM'=\bigvee\{\eK\in \fB_a;\;\eK\not\subseteq\eM\}=\bigvee\{\eK\in \fB_a;\;\eM\subseteq\eK'\}.$$
\end{corollary}

\begin{proof}
By the definition of atomic Boolean subspace lattices, we have
$\eM'=\bigvee\{\eK\in \fB_a;\;\eK\subseteq\eM'\}$. By \Cref{lem02}, an atom $\eK\in \fB_a$ is a
subspace of $\eM'$ if and only if it is not a subspace of $\eM$. Hence,
$\eM'=\bigvee\{\eK\in \fB_a;\;\eK\not\subseteq\eM\}$. Of course, by \Cref{lem02}, 
$\eK\not\subseteq\eM$ if and only if $\eM\subseteq\eK'$.
\end{proof}

The proof of the following theorem is an adaptation of Halmos' proof \cite[Theorem 1]{Hal}.

\begin{theorem} \label{theo02}
Every atomic Boolean subspace lattice is reflexive.
\end{theorem}

\begin{proof}
Let $\fB$ be an arbitrary atomic Boolean subspace lattice. Suppose that $\eM\in \Lat\Alg(\fB)$.
Let $\eM_0=\bigvee\{ \eK\in \fB_a;\; \eK\subseteq \eM\}$ and 
$\eM_1=\bigwedge\{\eK';\; \eK\in \fB_a,\;\eM\subseteq\eK'\}$. Of course, $\eM_0,\eM_1\in \fB$ and 
$\eM_0\subseteq \eM\subseteq\eM_1$. By \Cref{cor01},
$\eM_{0}^{\prime}=\bigvee\{ \eK\in \fB_a;\; \eM\subseteq \eK'\}$, and
by DeMorgan law,  $\eM_{1}^{\prime}=\bigvee\{ \eK\in \fB_a;\;\eM\subseteq\eK'\}$.
Thus, $\eM_{0}^{\prime}=\eM_{1}^{\prime}$ and therefore $\eM_0=\eM_1$ which 
gives $\eM=\eM_0\in \fB$.
\end{proof}

With some additional work, we can prove more. Recall that $\Rk1(\Alg(\fB))$ denotes the set of rank-one
operators in $\Alg(\fB)$.

\begin{theorem} \label{theo08}
Every atomic Boolean subspace lattice $\fB$ is strongly reflexive; in particular, $\fB=\Lat\bigl(\Rk1(\Alg(\fB))\bigr)$.
\end{theorem}

\begin{proof}
We have to prove that $\eN_{*}=\eN$, for all $\eN\in \fB$. Since $\eX_{*}=\eX$ and
$\{0\}_{*}=\bigwedge\limits_{\substack{\eM\in\fB\\ \eM\ne \{0\}}}\eM_{-}\subseteq
\bigwedge\limits_{\eK\in\fB_a}\eK_{-}=\{ 0\}$ we may suppose that $\eN$ is non-trivial.
Note that
\begin{equation*}
\begin{split}
\eN_{*}&=\bigwedge\{ \eM_{-};\;\eM\in \fB\;\;\text{and}\;\; \eM\not\subseteq \eN\}\\
&= \left(\bigwedge\{ \eM_{-};\;\eM\in \fB_a\;\;\text{and}\;\; \eM\not\subseteq \eN\}\right)\bigwedge
\left(\bigwedge\{ \eM_{-};\;\eM\in \fB\setminus\fB_a\;\;\text{and}\;\; \eM\not\subseteq \eN\}\right).
\end{split}
\end{equation*}
By \Cref{lem04}, $\eM_{-}=\eX$ if $\eM\ne \{0\}$ is not an atom. Hence, by \Cref{lem04} and \Cref{lem03},
$$\eN_{*}=\bigwedge\{ \eM';\; \eM\in \fB_a\;\;\text{and}\;\; \eM\subseteq \eN'\}. $$

Let $\eL\in \fB_a$ be such that $\eL\subseteq \eN$. If $\eM\in \fB_a$ and $\eM\subseteq \eN'$,
then $\eL\subseteq \eM'$. Hence,
$\eL\subseteq \bigwedge\{ \eM';\; \eM\in \fB_a\;\;\text{and}\;\; \eM\subseteq \eN'\}=\eN_{*}$.
On the other hand, if $\eL\in \fB_a$ and $\eL\subseteq \eN_{*}$, then $\eL\subseteq \eM'$, for all
$\eM\in \fB_a$ such that $\eM\subseteq \eN'$. It follows that $\eL\ne \eM$, for all
$\eM\in \fB_a$ such that $\eM\subseteq \eN'$. Hence, $\eL\not\subseteq \eN'$ and therefore $\eL\subseteq\eN$.
We conclude that $\eN_{*}=\eN$, that is, $\fB$ is strongly reflexive.
By \Cref{cor04}, strong reflexivity gives $\fB=\Lat\bigl(\Rk1(\Alg(\fB))\bigr)$.
\end{proof}

In \cite[Problem 9]{Hal1} Halmos raised a question if every Boolean subspace lattice is reflexive. 
It turns out that a Boolean subspace lattice that is not atomic does not need to be reflexive. Conway \cite{Con}
constructed a non-reflexive Boolean subspace lattice in the infinite-dimensional separable Hilbert space,
more precisely, in a subspace of the Hardy space $H^2$.
A brief summary of his technique (see \cite{Hal2}) is as follows. 
Let $\mu$ be a positive singular measure on the unit circle $\bT$. For each Borel subset $E \subseteq \bT$, 
let $\varphi_E$ be the inner function on the open unit disc $\bD$, defined by
$$\varphi_E(z)=\exp\left(-\int_E\frac{\lambda+z}{\lambda-z}d\mu(\lambda)\right).$$
If $\eM_E$ is the orthogonal complement of $\varphi_{\bT} H^2$ in $\varphi_E H^2$, that is,
$\eM_E=\varphi_E H^2 \cap (\varphi_{\bT} H^2)^\perp$, then $\eM_E$ is a subspace of 
$H^2 \cap (\varphi_{\bT} H^2)^\perp$. The correspondence $E\mapsto \eM_E$ is a Boolean anti-isomorphism
and its range $\fB\subseteq \fC\left( H^2 \cap (\varphi_{\bT} H^2)^\perp)\right)$ is a Boolean
subspace lattice. It turns out that $\fB$ is reflexive if and only if the measure $\mu$ is purely atomic.

\section{Ordinal sums} \label{Sec6}
\setcounter{theorem}{0}

Let $L_1$ and $L_2$ be abstract lattices (containing $0$ and $1$). The ordinal sum $L_1\osum L_2$ is the 
lattice whose elements are the elements of the disjoint union of $L_1$ and $L_2$ with the additional convention 
that the unit of  $L_1$ and the zero of $L_2$ are amalgamated. The ordering in $L_1\osum L_2$ is given as follows.
The relation between $x, y\in L_1$ or $x, y\in L_2$ is not changed. But for $x\in L_1$ and $y\in L_2$ we have
$x<y$, except when $x=1\in L_1$ and $y=0\in L_2$ (in this case $x=y$ as we amalgamated these two elements).
See \cite{Hal} for this definition. An example of the ordinal sum of two lattices is in \Cref{fig03}. 

\begin{center}
\begin{tikzpicture}[scale=0.8] 
\draw (0cm,1cm)--(0.8cm,1.5cm);
\draw (0.8cm,1.5cm)--(0.8cm,2.5cm);
\draw (0.8cm,2.5cm)--(0cm,3cm);
\draw (0cm,1cm)--(-0.8cm,1.5cm);
\draw (-0.8cm,1.5cm)--(-0.8cm,2.5cm);
\draw (-0.8cm,2.5cm)--(0cm,3cm);
\fill[white] (0cm,1cm) circle (2pt);
\draw (0cm,1cm) circle (2pt) node[below, scale=0.9] {$0$};
\fill[white] (-0.8cm,1.5cm) circle (2pt);
\draw (-0.8cm,1.5cm) circle (2pt) node[left, scale=0.9] {$a$};
\fill[white] (-0.8cm,2.5cm) circle (2pt);
\draw (-0.8cm,2.5cm) circle (2pt) node[left, scale=0.9] {$c$};
\fill[white] (0.8cm,1.5cm) circle (2pt);
\draw (0.8cm,1.5cm) circle (2pt) node[below right, scale=0.9] {$b$};
\fill[white] (0.8cm,2.5cm) circle (2pt);
\draw (0.8cm,2.5 cm) circle (2pt) node[above right, scale=0.9] {$d$};
\fill[white] (0cm,3cm) circle (2pt);
\draw (0cm,3cm) circle (2pt) node[above, scale=0.9] {$1$};
\draw (0cm,-0.8cm) node[below, scale=0.9] {$L_1$};
\draw (1.3cm,2cm)  node {$\osum$};
\draw (3cm,1cm)--(4cm,2cm);
\draw (4cm,2cm)--(3cm,3cm);
\draw (3cm,1cm)--(2cm,2cm);
\draw (2cm,2cm)--(3cm,3cm);
\fill[white] (3cm,1cm) circle (2pt);
\draw (3cm,1cm) circle (2pt) node[below, scale=0.9] {$0$};
\fill[white] (2cm,2cm) circle (2pt);
\draw (2cm,2cm) circle (2pt) node[left, scale=0.9] {$e$};
\fill[white] (4cm,2cm) circle (2pt);
\draw (4cm,2cm) circle (2pt) node[right, scale=0.9] {$f$};
\fill[white] (3cm,3cm) circle (2pt);
\draw (3cm,3cm) circle (2pt) node[above, scale=0.9] {$1$};
\draw (3cm,-0.8cm) node[below, scale=0.9] {$L_2$};
\draw (4.9cm,2cm)  node {$=$};
\draw (6.5cm,0cm)--(7.3cm,0.5cm);
\draw (7.3cm,0.5cm)--(7.3cm,1.5cm);
\draw (7.3cm,1.5cm)--(6.5cm,2cm);
\draw (6.5cm,0cm)--(5.7cm,0.5cm);
\draw (5.7cm,0.5cm)--(5.7cm,1.5cm);
\draw (5.7cm,1.5cm)--(6.5cm,2cm);
\fill[white] (6.5cm,0 cm) circle (2pt);
\draw (6.5cm,0 cm) circle (2pt) node[below, scale=0.9] {$0$};
\fill[white] (5.7cm,0.5cm) circle (2pt);
\draw (5.7cm,0.5cm) circle (2pt) node[left, scale=0.9] {$a$};
\fill[white] (5.7cm,1.5cm) circle (2pt);
\draw (5.7cm,1.5cm) circle (2pt) node[left, scale=0.9] {$c$};
\fill[white] (7.3cm,0.5cm) circle (2pt);
\draw (7.3cm,0.5cm) circle (2pt) node[below right, scale=0.9] {$b$};
\fill[white] (7.3cm,1.5cm) circle (2pt);
\draw (7.3cm,1.5 cm) circle (2pt) node[above right, scale=0.9] {$d$};
\fill[white] (6.5cm,2cm) circle (2pt);
\draw (6.5cm,1.95cm) node[below, scale=0.9] {$p$};
\fill[black] (6.5cm,2cm) circle (2pt);
\draw (6.5cm,2cm)--(7.5cm,3cm)--(6.5cm,4cm)--(5.5cm,3cm)--cycle;
\fill[white] (5.5cm,3cm) circle (2pt);
\draw (5.5cm,3cm) circle (2pt) node[left, scale=0.9] {$e$};
\fill[white] (7.5cm,3cm) circle (2pt);
\draw (7.5cm,3cm) circle (2pt) node[right, scale=0.9] {$f$};
\fill[white] (6.5cm,4cm) circle (2pt);
\draw (6.5cm,4cm) circle (2pt) node[above, scale=0.9] {$1$};
\draw (6.5cm,-0.8cm) node[below, scale=0.9] {$L_1\osum L_2$};
\end{tikzpicture}

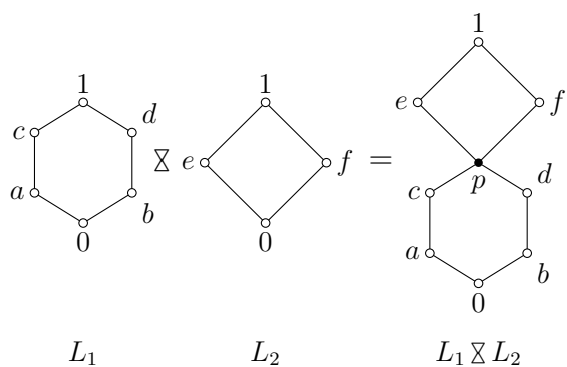
\captionof{figure}{Ordinal sum of two lattices.}   \label{fig03}
\end{center}

\noindent
As we see, the Hasse diagram of the ordinal sum of lattices $L_1$ and $L_2$ in \Cref{fig03} is obtained by pasting the Hasse
diagrams of $L_1$ and $L_2$, equating the greatest element of $L_1$ with the least element of $L_2$.

The notion of a pinch point has sense in an abstract lattice $L$, too. We say that an element $p\in L$ is 
a pinch point if $a\leq p$ or $p\leq a$, for every $a\in L$. 
For instance, in \Cref{fig03}, element $p$ in the ordinal sum $L_1\osum L_2$ is a pinch point.

Let $\fL_i$ be a realization of a lattice $L_i$ as a subspace lattice in a Banach space $\eX_i$ 
and let $\phi_i\colon L_i\to \fL_i$ be the corresponding isomorphism $(i=1,2)$. Then the ordinal sum 
$L_1\osum L_2$ is realizable in the direct sum $\eX_1\oplus\eX_2$. We can realize it as follows. 
We define $\phi\colon L_1\osum L_2\to \fC(\eX_1\oplus \eX_2)$ by
$$ \phi(a)=\left\{ \begin{array}{cl}
\phi_1(a)\oplus \{0\},\quad & a\in L_1,\\
\eX_1\oplus \phi_2(a),\quad & a\in L_2.
\end{array} \right. $$
It is obvious from the definition that $\phi$ is an isomorphism from $L_1\osum L_2$ onto a subspace lattice in
$\fC(\eX_1\oplus \eX_2)$. We will write $\fL_1\osum\fL_2$ for this realization of the ordinal sum $L_1\osum L_2$.

\begin{lemma} \label{lem12}
Operator $A\in \pB(\eX_1\oplus \eX_2)$ is in $\Alg(\fL_1\osum\fL_2)$ if and only if there exist
$A_{1}\in \Alg(\fL_1)$, $A_2\in \Alg(\fL_2)$, and $B\in \pB(\eX_2,\eX_1)$ such that
\begin{equation} \label{eq03}
A(x_1\oplus x_2)=(A_1 x_1+Bx_2)\oplus A_2x_2,\qquad \text{for all}\; x_1\oplus x_2\in \eX_1\oplus\eX_2.
\end{equation}
\end{lemma}

\begin{proof}
It is easy to check that for arbitrary operators $A_{1}\in \Alg(\fL_1)$, $A_2\in \Alg(\fL_2)$, and 
$B\in \pB(\eX_2,\eX_1)$ the operator $A$, which is given by 
$A(x_1\oplus x_2)=(A_1 x_1+Bx_2)\oplus A_2x_2$, for all $x_1\oplus x_2\in \eX_1\oplus\eX_2$,
belongs to $\Alg(\fL_1\osum\fL_2)$. For the opposite implication, note that an arbitrary operator
$A\in \pB(\eX_1\oplus\eX_2)$ is of the form 
$A(x_1\oplus x_2)=(A_1 x_1+Bx_2)\oplus (Cx_1+A_2x_2)$, for all $x_1\oplus x_2\in \eX_1\oplus\eX_2$,
where $A_1\in \pB(\eX_1)$, $A_2\in \pB(\eX_2)$, $B\in \pB(\eX_2,\eX_1)$, and $C\in \pB(\eX_1,\eX_2)$
are arbitrary. If $A\in \Alg(\fL_1\osum\fL_2)$, then it follows from $A(\eX_1\oplus\{ 0\})\subseteq 
\eX_1\oplus\{ 0\}$ that $C=0$. Let $\eM_1\in \fL_1$ be arbitrary. Then $\eM_1\oplus\{ 0\}\in \fL_1\osum\fL_2$
and therefore $A(x_1\oplus 0)=A_1 x_1\oplus 0\in \eM_1\oplus\{ 0\}$, for every $x_1\in \eM_1$. It follows that
$A_1\eM_1\subseteq \eM_1$ and consequently $A_1\in \Alg(\fL_1)$. Similarly, if $\eM_2\in \fL_2$,
then $\eX_1\oplus\eM_2\in \fL_1\osum\fL_2$. Hence, for all $x_1\in \eX_1$ and $x_2\in \eM_2$ we
have $A(x_1\oplus x_2)=(A_1x_1+Bx_2)\oplus A_2 x_2\in \eX_1\oplus\eM_2$. This gives 
$A_2\eM_2\subseteq \eM_2$ and therefore $A_2\in \Alg(\fL_2)$.
\end{proof}

As the following theorem shows, the ordinal sum of two reflexive subspace lattices is reflexive. Our proof
is based on the proof of \cite[Theorem 2]{Hal}.

\begin{theorem} \label{theo05}
Let $\fL_1\subseteq \fC(\eX_1)$ and $\fL_2\subseteq \fC(\eX_2)$ be subspace lattices. The ordinal sum 
$\fL_1\osum\fL_2\subseteq \fC(\eX_1\oplus\eX_2)$ is reflexive if and only if $\fL_1$ and $\fL_2$ are reflexive.
\end{theorem}

\begin{proof}
Suppose that $\fL_1$ and $\fL_2$ are reflexive and let $\eM\in \Lat\Alg(\fL_1\osum\fL_2)$ be arbitrary.
If $\eM=\eM_1\oplus\{ 0\}$, for some $\eM_1\in \fC(\eX_1)$, then
$A(x_1\oplus 0)=A_1 x_1\oplus 0\in \eM_1\oplus\{ 0\}$, for all $x_1\in \eM_1$ and every operator
$A\in\Alg(\fL_1\osum\fL_2)$, which is given by \eqref{eq03}. Hence, $A_1\eM_1\subseteq\eM_1$,
for all $A_1\in \Alg(\fL_1)$, by \Cref{lem12}. Since $\fL_1$ is reflexive we have $\eM_1\in \fL_1$ and
therefore $\eM\in \fL_1\osum\fL_2$. Assume now that there exists $x_1\oplus x_2\in \eM$ such that $x_2\ne 0$.
By \Cref{lem12}, $Bx_2\oplus 0\in \eM$, for all $B\in \pB(\eX_2,\eX_1)$. It follows that $\eX_1\oplus\{ 0\}
\subseteq \eM$. Hence, $\eM=\eX_1\oplus \eM_2$, where $\eM_2\in \fC(\eX_2)$. Let $y_2\in \eM_2$ be
arbitrary. Then $0\oplus y_2\in \eM$ and $0\oplus A_2y_2\in \eM$, for all $A_2\in \Alg(\fL_2)$, by \Cref{lem12}.
Hence, $A_2\eM_2\subseteq \eM_2$, for all $A_2\in \Alg(\fL_2)$. Since $\fL_2$ is reflexive we have
$\eM_2\in \fL_2$ and consequently $\eM\in \fL_1\osum\fL_2$.

To prove the opposite implication, assume that $\fL_1\osum\fL_2$ is reflexive. For $\eM_1\in \Lat\Alg(\fL_1)$,
let $\eM=\eM_1\oplus\{ 0\}$. Let $A\in \Alg(\fL_1\osum\fL_2)$ be given by \eqref{eq03}. Then
$A(x_1\oplus 0)=A_1x_1\oplus 0\in \eM_1\oplus \{ 0\}$, for all $x_1\in \eM_1$. Hence, 
$\eM\in\Lat\Alg(\fL_1\osum\fL_2)=\fL_1\osum\fL_2$ and therefore $\eM_1\in \fL_1$. If $\eM_2\in \Lat\Alg(\fL_2)$,
then let $\eN=\eX_1\oplus \eM_2$. Again, if $A\in \Alg(\fL_1\osum\fL_2)$ is given by \eqref{eq03}, then
$A(x_1\oplus x_2)=(A_1x_1+Bx_2)\oplus A_2x_2\in \eN$, for all $x_1\in x_2\in \eN$. Hence, $\eN\in 
\fL_1\osum\fL_2$ and therefore $\eM_2\in \fL_2$.
\end{proof}

It is clear that the ordinal sum of subspace lattices $\fL_1$ and $\fL_2$ has a pinch point at the position 
where $\fL_1$ and $\fL_2$ are amalgamated. On the other hand, a subspace lattice
$\fL$ with a pinch point $\eP$ is an ordinal sum of two subspace lattices that are amalgamated at $\eP$.
To see this, let $\fL\subseteq \fC(\eX)$ be a subspace lattice and suppose that $\eP\in \fL$ is a pinch point. 
To avoid trivialities, assume that $\{0\}\ne\eP\ne\eX$. If $\eM\in \fL$, then either $\eM\subseteq \eP$ or 
$\eP\subseteq \eM$.
Denote by $\fL_{\eP}$ the family of all those subspaces in $\fL$ which are contained in $\eP$. It is clear
that $\fL_{\eP}$ can be considered as a subset of $\fC(\eP)$. If $\eP\subseteq \eM$, then let
$\langle\eM\rangle=\{ \langle x\rangle;\; x\in \eM\}\subseteq \eX/\eP$, where 
$\langle x\rangle=\{ x+p;\; p\in \eP\}$ is the coset in $\eX/\eP$
determined by $x$. Of course,  $\langle\eM\rangle$ is a subspace of the quotient space $\eX/\eP$. Let 
$\fL^{\eP}=\{ \langle \eM\rangle;\; \eM\in \fL\;\text{and}\; \eP\subseteq \eM\}$.

\begin{lemma} \label{lem09}
$\fL_{\eP}$ is a subspace lattice in $\fC(\eP)$ and $\fL^{\eP}$ is a subspace lattice in $\fC(\eX/\eP)$.
\end{lemma}

\begin{proof}
It is obvious that $\fL_{\eP}$ is a subspace lattice in $\fC(\eP)$. To see that $\fL^{\eP}$ is a subspace
lattice in $\fC(\eX/\eP)$, note first that $\langle\eP\rangle$ is equal to $\{ 0\}$ in $\eX/\eP$ and $\langle\eX\rangle$
is equal to $\eX/\eP$.
Hence, $\{ 0\}$ and $\eX/\eP$ are in $\fL^{\eP}$. Let $\{ \langle\eM_j\rangle;\; j\in J\}$ be an arbitrary family of
subspaces in $\fL^{\eP}$. We want to show that 
$\bigwedge\limits_{j\in J}\langle\eM_j\rangle= \langle\bigwedge\limits_{j\in J}\eM_j\rangle$ and 
$\bigvee\limits_{j\in J}\langle\eM_j\rangle= \langle\bigvee\limits_{j\in J}\eM_j\rangle$ since this implies that
$\bigwedge\limits_{j\in J}\langle\eM_j\rangle$ and $\bigvee\limits_{j\in J}\langle\eM_j\rangle$ are in $\fL^{\eP}$.

Suppose that $x\in \eX$ is such that $\langle x\rangle\in \bigwedge\limits_{j\in J}\langle\eM_j\rangle$. 
Then, for every $j\in J$, there
exists $x_j\in \eM_j$ such that $\langle x\rangle=\langle x_j\rangle$, that is $\{ x+p;\; p\in \eP\}=
\{ x_j+p;\; p\in \eP\}\subseteq \eM_j$
since $\eP\subseteq \eM_j$. Hence, $x\in \eM_j$, for all $j\in J$, and therefore 
$\langle x\rangle\subseteq \bigwedge\limits_{j\in J}\eM_j$ which gives 
$\langle x\rangle\in \langle\bigwedge\limits_{j\in J}\eM_j\rangle$.
On the other hand, if $\langle x\rangle\in \langle\bigwedge\limits_{j\in J}\eM_j\rangle$, then 
$\langle x\rangle\subseteq \bigwedge\limits_{j\in J}\eM_j+\eP=\bigwedge\limits_{j\in J}\eM_j$ and consequently
$x\in \eM_j$, for all $j\in J$. Hence $\langle x\rangle\subseteq \eM_j=\eM_j+\eP$ and therefore 
$\langle x\rangle\in \langle\eM_j\rangle$, for all $j\in J$.
Thus, $\langle x\rangle\in \bigwedge\limits_{j\in J}\langle \eM_j\rangle$.

Since $\bigvee\limits_{j\in J}\langle \eM_j\rangle\supseteq \langle\eM_i\rangle=\eM_i+\eP$, for every $i\in J$,
we have $\bigvee\limits_{j\in J}\langle\eM_j\rangle\supseteq (\eM_{i_1}+\eP)+\cdots+(\eM_{i_k}+\eP)=
(\eM_{i_1}+\cdots+\eM_{i_k})+\eP$, for every finite family $\eM_{i_1},\ldots,\eM_{i_k}$ of subspaces in the family 
$\{ \eM_j;\; j\in J\}$. It follows that 
$\bigvee\limits_{j\in J}\langle \eM_j\rangle\supseteq \bigvee\limits_{j\in J}\eM_j+\eP=
\langle\bigvee\limits_{j\in J}\eM_j\rangle$.
On the other hand, $\langle\eM_i\rangle\subseteq \langle\bigvee\limits_{j\in J}\eM_j\rangle$, for every 
$i\in J$, and therefore
$\bigvee\limits_{i\in J}\langle\eM_i\rangle\subseteq \langle\bigvee\limits_{j\in J}\eM_j\rangle$. 
Thus, $\bigvee\limits_{j\in J}\langle\eM_j\rangle= \langle\bigvee\limits_{j\in J}\eM_j\rangle\in \fL^{\eP}$.
\end{proof}

\begin{proposition} \label{prop02}
Let $\fL\subseteq \fC(\eX)$ be a subspace lattice and let $\eP\in \fL$ be a pinch point such
that $\{0\}\ne\eP\ne\eX$. Then $\fL$ is isomorphic to $\fL_{\eP}\osum\fL^{\eP}$.
\end{proposition}

\begin{proof}
Note first that $\fL_{\eP}\osum\fL^{\eP}$ is a subspace lattice in $\fC\bigl(\eP\oplus (\eX/\eP)\bigr)$.
Define $\phi\colon\fL\to \fL_{\eP}\osum\fL^{\eP}$ as follows. Let $\eM\in\fL$ be arbitrary.
Since $\eP$ is a pinch point we have either $\eM\subseteq \eP$ or $\eP\subsetneq \eM$.
If $\eM\subseteq \eP$, then let $\phi(\eM)=\eM\oplus \langle\{ 0\}\rangle$, and if $\eP\subsetneq \eM$,
then let $\phi(\eM)=\eP\oplus \langle\eM\rangle$. It is obvious that $\phi$ is a bijection such that
$\phi(\eM_1)\subseteq \phi(\eM_2)$ if and only if $\eM_1\subseteq \eM_2$.
\end{proof}

\section{Small lattices} \label{Sec7}
\setcounter{theorem}{0}

In this section, we consider lattices with a few elements. Since we assume that every lattice has $0$ and $1$
all lattices have at least two elements and that with two elements is trivial. There is only one lattice with three
elements and it is a nest. Of course, it is realizable in every Banach space of dimension at least $2$ and every 
realization of it is strongly reflexive.

There are two non-isomorphic lattices with four elements: the chain (nest) and the diamond, see \Cref{fig04}.

\begin{center}
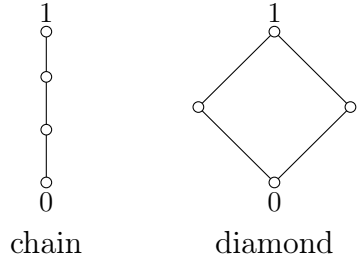

\begin{tikzpicture}[scale=1] 
\draw (-3cm,0cm)--(-3cm,2cm);
\fill[white] (-3cm,0cm) circle (2pt);
\draw (-3cm,0 cm) circle (2pt) node[below, scale=0.9] {$0$};
\fill[white] (-3cm,0.7cm) circle (2pt);
\draw (-3cm,0.7cm) circle (2pt);
\fill[white] (-3cm,1.4cm) circle (2pt);
\draw (-3cm,1.4cm) circle (2pt);
\fill[white] (-3cm,2cm) circle (2pt);
\draw (-3cm,2cm) circle (2pt) node[above, scale=0.9] {$1$};
\draw (-3cm,-0.5cm) node[below] {chain};
\draw (0cm,0cm)--(1cm,1cm);
\draw (1cm,1cm)--(0cm,2cm);
\draw (0cm,0cm)--(-1cm,1cm);
\draw (-1cm,1cm)--(0cm,2cm);
\fill[white] (0cm,0 cm) circle (2pt);
\draw (0cm,0 cm) circle (2pt) node[below, scale=0.9] {$0$};
\fill[white] (-1cm,1 cm) circle (2pt);
\draw (-1cm,1cm) circle (2pt);
\fill[white] (1cm,1cm) circle (2pt);
\draw (1cm,1cm) circle (2pt);
\fill[white] (0cm,2cm) circle (2pt);
\draw (0cm,2 cm) circle (2pt) node[above, scale=0.9] {$1$};
\draw (0cm,-0.5cm) node[below] {diamond};
\end{tikzpicture}
\captionof{figure}{Lattices with four elements.}   \label{fig04}
\end{center}

The chain is realizable in every Banach space of dimension at least $3$ and its realization is strongly reflexive. The diamond can be realized as a subspace lattice in every 
Banach space of dimension at least $2$. Since a realization of the diamond is atomic and Boolean it is strongly reflexive, too.

There are five isomorphism types of lattices with exactly five elements: the chain (nest),
the pendulum, the pendulum upside down, the pentagon, and the double triangle, see Figures \ref{fig05}, \ref{fig01},
and \ref{fig06}.

\begin{center}
\begin{tikzpicture}[scale=1] 
\draw (-0.5cm,0cm)--(-0.5cm,2.8cm);
\fill[white] (-0.5cm,0cm) circle (2pt);
\draw (-0.5cm,0 cm) circle (2pt) node[below, scale=0.9] {$0$};
\fill[white] (-0.5cm,0.7cm) circle (2pt);
\draw (-0.5cm,0.7cm) circle (2pt);
\fill[white] (-0.5cm,1.4cm) circle (2pt);
\draw (-0.5cm,1.4cm) circle (2pt);
\fill[white] (-0.5cm,2.1cm) circle (2pt);
\draw (-0.5cm,2.1cm) circle (2pt);
\fill[white] (-0.5cm,2.8cm) circle (2pt);
\draw (-0.5cm,2.8cm) circle (2pt) node[above, scale=0.9] {$1$};
\draw (-0.5cm,-0.5cm) node[below] {chain};
\draw (3cm,0cm)--(2cm,1cm)--(3cm,2cm)--(4cm,1cm)--cycle;
\draw (3cm,2.8cm)--(3cm,2cm);
\fill[white] (3cm,0 cm) circle (2pt);
\draw (3cm,0 cm) circle (2pt) node[below, scale=0.9] {$0$};
\fill[white] (2cm,1cm) circle (2pt);
\draw (2cm,1cm) circle (2pt);
\fill[white] (4cm,1cm) circle (2pt);
\draw (4cm,1cm) circle (2pt);
\fill[white] (3cm,2cm) circle (2pt);
\draw (3cm,2cm) circle (2pt);
\fill[white] (3cm,2.8cm) circle (2pt);
\draw (3cm,2.8cm) circle (2pt) node[above, scale=0.9] {$1$};
\draw (3cm,-0.5cm) node[below] {pendulum};
\draw (7cm,0.8cm)--(6cm,1.8cm)--(7cm,2.8cm)--(8cm,1.8cm)--cycle;
\draw (7cm,0cm)--(7cm,0.8cm);
\fill[white] (7cm,0 cm) circle (2pt);
\draw (7cm,0 cm) circle (2pt) node[below, scale=0.9] {$0$};
\fill[white] (7cm,0.8cm) circle (2pt);
\draw (7cm,0.8cm) circle (2pt);
\fill[white] (6cm,1.8cm) circle (2pt);
\draw (6cm,1.8cm) circle (2pt);
\fill[white] (8cm,1.8cm) circle (2pt);
\draw (8cm,1.8cm) circle (2pt);
\fill[white] (7cm,2.8cm) circle (2pt);
\draw (7cm,2.8cm) circle (2pt) node[above, scale=0.9] {$1$};
\draw (7cm,-0.5cm) node[below] {pendulum};
\draw (7cm,-0.9cm) node[below] {upside down};
\end{tikzpicture}

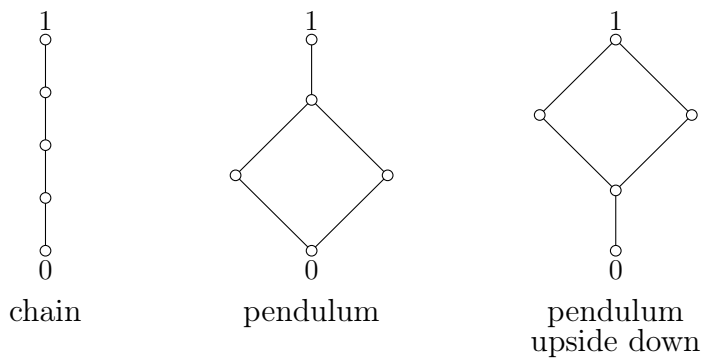
\captionof{figure}{Three five-element lattices.}   \label{fig05}
\end{center}

Again, the chain with five elements is realizable in every Banach space of dimension $4$ or more
and every realization of it is strongly reflexive.
The pendulums can be realized as subspace lattices in a Banach space of dimension at least $3$. 
Since the pendulums are ordinal sums of the trivial 
chain and the diamond a realization of a pendulum is reflexive, by \Cref{theo05}. In what follows
we will consider realizations of the pentagon and the double triangle.

\subsection{Pentagon}
The Hasse diagram of this abstract lattice $P=\{0,a,b,c,1\}$ is in \Cref{fig01}.

\begin{center}
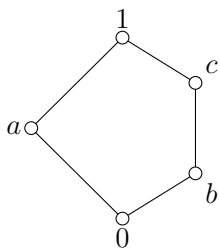

\begin{tikzpicture}[scale=1.2] 
\draw (0cm,0cm)--(0.8cm,0.5cm);
\draw (0.8cm,0.5cm)--(0.8cm,1.5cm);
\draw (0.8cm,1.5cm)--(0cm,2cm);
\draw (0cm,0cm)--(-1cm,1cm);
\draw (-1cm,1cm)--(0cm,2cm);
\fill[white] (0cm,0 cm) circle (2pt);
\draw (0cm,0 cm) circle (2pt) node[below, scale=0.9] {$0$};
\fill[white] (-1cm,1 cm) circle (2pt);
\draw (-1cm,1cm) circle (2pt) node[left, scale=0.9] {$a$};
\fill[white] (0.8cm,0.5cm) circle (2pt);
\draw (0.8cm,0.5cm) circle (2pt) node[below right, scale=0.9] {$b$};
\fill[white] (0.8cm,1.5cm) circle (2pt);
\draw (0.8cm,1.5 cm) circle (2pt) node[above right, scale=0.9] {$c$};
\fill[white] (0cm,2cm) circle (2pt);
\draw (0cm,2 cm) circle (2pt) node[above, scale=0.9] {$1$};
\end{tikzpicture}
\captionof{figure}{The pentagon lattice.}   \label{fig01}
\end{center}
As we can see, $P$ is a complemented lattice however the complement of $a$ is not unique. 
On the other hand, the pentagon is not modular. Indeed, we have $b\leq c$, but
$b\vee(a\wedge c)=b$ and $(b\vee a)\wedge c=c$, that is, $b\vee(a\wedge c)\ne (b\vee a)\wedge c$.
By \Cref{lem06}, every subspace lattice in a finite-dimensional Banach space is modular. 
Since the pentagon is not modular it is not realizable in a finite-dimensional Banach space.
However, it can be realized as a subspace lattice in some infinite-dimensional complex Banach space $\eX$.

Let $\fP=\bigl\{ \{ 0\},\eM,\eN,\widehat{\eN},\eX\bigr\}$, where $\eN\subsetneq \widehat{\eN}$, 
$\eM\vee\eN=\eX$, and $\eM\wedge\widehat{\eN}=\{0\}$, be a realization of the pentagon as 
a subspace lattice in an infinite-dimensional Banach space $\eX$. The {\em gap-dimension} of $\fP$ is 
the dimension of the quotient space $\widehat{\eN}/\eN$.
Following Halmos \cite{Hal}, we will show that for every $k\in \bN$ there exists a realization $\fP_k$ of the
pentagon with the gap-dimension $k$ whenever there exist a quasi-complemented pair of subspaces $(\eM,\eN)$ in 
$\fC(\eX)$ such that $\eM+\eN\ne \eX$, that is, $\eN$ is not a complement of $\eM$. 
It is well-known that subspaces with these
properties exist in every infinite-dimensional Hilbert space. On the other hand, Murray \cite{Mur} and Mackey
proved that every subspace a separable Banach space $\eX$ is quasi-complemented. Hence, if $\eX$ is not 
isomorphic to a Hilbert space,
then there exist subspaces $\eM, \eN\in\fC(\eX)$ such that $ \eM\wedge\eN=\{ 0\}$, $\eM\vee\eN=\eX$, but 
$\eM+\eN\ne \eX$, by the famous result of Lindenstrauss and Tzafriri \cite{LT}. 

Let $\eX$ be a Banach space such that there exists a quasi-complemented pair of subspaces $(\eM, \eN)$ and 
$\eM+\eN\ne \eX$.
It follows that $\eM$ and $\eN$ are proper infinite-dimensional subspaces. Let $u_1\in \eX$ be such that 
$u_1\not\in \eM+\eN$ and let $[u_1]$ be the one-dimensional subspace spanned by $u_1$.
Denote $\eN_1=\eN\vee [u_1]$. It is clear that $\eN\cap [u_1]=\{ 0\}$ and therefore
$\eN_1=\eN\oplus [u_1]$. Of course, $\eM\vee \eN_1=\eX$. We claim that $\eM\wedge \eN_1=\{ 0\}$.
Indeed, suppose the opposite, that is, there is $0\ne x\in \eM\wedge\eN_1$. Since $\eM\wedge\eN=\{ 0\}$ we
see that $x=y+\alpha u_1$, where $y\in \eN$ and $\alpha\ne 0$. Hence $u_1=\alpha ^{-1}(x-y)$. But 
$\alpha ^{-1}(x-y)\in \eM+\eN$ and therefore $u_1\in \eM+\eN$ which is a contradiction. We conclude that
$\eM\wedge\eN_1=\{ 0\}$. This shows that $\fP_1=\bigl\{ \{ 0\}, \eM,\eN,\eN_1,\eX\bigr\}$ is a subspace 
lattice in $\eX$; it is a realization of the pentagon and its gap-dimension is $1$.

We have seen that $\eM\wedge\eN_1=\{ 0\}$ and $\eM\vee\eN_1=\eX$. Let us show that $\eM+\eN_1\ne\eX$.
Indeed, if $\eM+\eN_1=\eX$, then $\eN_1$ can be
identified with the quotient space $\eX/\eM$. Since $\eN_1=\eN\oplus [u_1]$, there exists a
functional $\theta\in (\eX/\eM)^*$ such that $\theta\in \eN^\perp$ and $\langle u_1,\theta\rangle=1$. The space
$(\eX/\eM)^*$ can be identified with $\eM^\perp$. Thus, there exists a functional $\vartheta\in \eX^*$
which annihilates $\eM$ and $\eN$ but $\langle u_1,\vartheta\rangle=1$. This is impossible since $\eM+\eN$ is
dense in $\eX$.
Subspace $\eN_1$ has the same properties as $\eN$. Thus, similarly, as before, we can construct a subspace
$\eN_2\in \fC(\eX)$ such that $ \eM\wedge\eN_2=\{ 0\}$, $\eM\vee\eN_2=\eX$, and $\eM+\eN_2\ne \eX$.
The codimension of $\eN$ in $\eN_2$ is $2$. Inductively we can get a subspace 
$\eN_k\in \fC(\eX)$ such that $\eN$ is a subspace in $\eN_k$ with codimension $k$. It is clear that
$\fP_k=\bigl\{ \{ 0\}, \eM,\eN,\eN_k,\eX\bigr\}$ is a realization of the pentagon with the gap-dimension $k$.

If $\fP=\bigl\{ \{ 0\},\eM,\eN,\widehat{\eN},\eX\bigr\}$ is a pentagon subspace lattice, then, for every 
$n\in \bN\cup\{\infty\}$, the subspace lattice 
$\fP^{(n)}=\bigl\{ \{ 0\}^{(n)},\eM^{(n)},\eN^{(n)},\widehat{\eN}^{(n)},$ $\eX^{(n)}\bigr\}$
is a realization of the pentagon, as well, by \Cref{lem13}. It is obvious that the gap-dimension of 
$\fP^{(n)}$ is $n$ times the gap-dimension of $\fP$. Thus, there exist pentagon subspace lattices of an
arbitrary gap-dimension.

\begin{theorem} \label{theo03}
A pentagon subspace lattice $\fP=\bigl\{ \{ 0\},\eM,\eN,\widehat{\eN},\eX\bigr\}$ is not strongly reflexive. 
Moreover, $\Lat\bigl(\Rk1(\Alg(\fP))\bigr)=\fP\cup[\eN,\widehat{\eN}]$.
\end{theorem}

\begin{proof}
Recall that $\Rk1(\Alg(\fP))$ denotes the subset of rank-one operators in $\Alg(\fP)$.
We have $\eN\subsetneq \widehat{\eN}$, $\eM\vee\eN=\eX$, and $\eM\wedge\widehat{\eN}=\{0\}$.
By \eqref{eq01}, $\eM_{-}=\{ 0\}\vee\eN\vee\widehat{\eN}=\widehat{\eN}$. Similarly, 
$\eN_{-}=\eM$ and $\widehat{\eN}_{-}=\eX$. It follows, see \eqref{eq04}, that $\eN_{*}=\widehat{\eN}\ne \eN$.
Hence, $\fP$ is not strongly reflexive.

By \Cref{lem01}, a rank-one operator $e\otimes\eta$ is in $\Alg(\fP)$ if and only if either $e\in\eM$ and
$\eta\in\widehat{\eN}^\perp$ or $e\in \eN$ and $\eta\in\eM^\perp$. In any case, it is not hard to see
that $\eK\in[\eN,\widehat{\eN}]$ is invariant for $e\otimes\eta$. Hence, $\fP\cup[\eN,\widehat{\eN}]
\subseteq \Lat\bigl(\Rk1(\Alg(\fP))\bigr)$.

Suppose now that $\eK\in \Lat\bigl(\Rk1(\Alg(\fP))\bigr)$. There is no loss of generality if we assume that
$\{0\}\ne \eK\ne\eX$. First, we will prove the following:
\begin{itemize}
\item[(\dag)]\quad if there exists $0\ne f\in \eK$ such that $f\not\in \widehat{\eN}$, then $\eM\subseteq \eK$.
\end{itemize}
Indeed, since $f\not\in \widehat{\eN}$ there exists $\eta\in \widehat{\eN}^\perp=(\eM_{-})^\perp$ such that 
$\langle f,\eta\rangle=1$. For every $0\ne e\in\eM$, the rank-one operator $e\otimes \eta$ is in $\Alg(\fP)$.
Hence, $(e\otimes \eta)f=e\in \eK$, for all $e\in \eM$, and therefore $\eM\subseteq \eK$. 

Now we consider two possibilities: (a) $\eK\subseteq \eM$ and 
(b) $\eK\not\subseteq \eM$. In the former case, we have $\eK\wedge \widehat{\eN}=\{ 0\}$ since 
$\eM\wedge \widehat{\eN}=\{ 0\}$. Hence, there exists $0\ne f\in \eK$ which is not in $\widehat{\eN}$ 
and therefore $\eM\subseteq \eK$, by (\dag). Thus, $\eK=\eM$ in this case.
Assume (b). Then there exist $g\in \eK\setminus \eM$ and  $\xi\in \eM^\perp=(\eN_{-})^\perp$
such that $\langle g,\xi\rangle=1$. Since, by \Cref{lem01}, $h\otimes \xi\in \Alg(\fP)$, for all $h\in \eN$,
we have $(h\otimes \xi)g=h\in \eK$, for all $h\in \eN$. This gives $\eN\subseteq \eK$. If there
existed $f\in \eK\setminus \widehat{\eN}$, then we would have $\eM\subseteq \eK$, by (\dag), and 
consequently $\eX=\eM\vee\eN\subseteq \eK\ne \eX$, which is impossible. Hence, we may conclude that
$\eK\subseteq \widehat{\eN}$.
\end{proof}

The following corollary is restated from \cite[Proposition 2.2]{KLL}.

\begin{corollary} \label{cor02}
If $\fP=\bigl\{ \{ 0\},\eM,\eN,\widehat{\eN},\eX\bigr\}$ is a pentagon subspace lattice with the gap-dimension $1$,
then $\fP^{(n)}=\bigl\{ \{ 0\}^{(n)},\eM^{(n)},\eN^{(n)},\widehat{\eN}^{(n)},\eX^{(n)}\bigr\}$
is reflexive, for every $n\in \bN\cup\{\infty\}$.
\end{corollary}

\begin{proof}
If the gap-dimension of $\fP$ is $1$, then $[\eN,\widehat{\eN}]=\{\eN,\widehat{\eN}\}$ and therefore
$\fP\cup[\eN,\widehat{\eN}]=\fP$ which gives 
$\fP\subseteq \Lat\Alg(\fP)\subseteq \Lat\bigl(\Rk1(\Alg(\fP))\bigr)=\fP$, by \Cref{theo03}.
Hence, $\fP$ is reflexive.

Let $n\in \bN\cup\{ \infty\}$, $n\geq 2$, be arbitrary and let $\eK\in \Lat\Alg\bigl(\fP^{(n)}\bigr)$.
For every integer $1\leq k\leq n$, let $P_k\colon \eX^{(n)}\to\eX^{(n)}$
be the projection onto the $k$th coordinate, that is, $P_k$ maps $(x_j)_{j=1}^{n}\in \eX^{(n)}$ into
$(y_j)_{j=1}^{n}\in \eX^{(n)}$, where $y_k=x_k$ and $y_j=0$ if $j\ne k$. It is clear that 
$P_k\in\Alg\bigl(\fP^{(n)}\bigr)$, for all integers $1\leq k\leq n$. Hence, $P_k\eK\subseteq \eK$. We may 
identify $P_k\eK$ with a subspace $\eK_k\subseteq \eX$. Then 
$\eK$ is a direct sum of subspaces $\eK_k\in \fC(\eX)$, $1\leq k\leq n$. 
For $A\in \pB(\eX)$ and an integer $1\leq k\leq n$, let operator $\widetilde{A}_k\in \pB\bigl(\eX^{(n)}\bigr)$
be given by $\widetilde{A}_k(x_j)_{j=1}^{n}=(y_j)_{j=1}^{n}$, where $y_k=Ax_k$ and $y_j=0$ if $j\ne k$.
It is obvious that $\widetilde{A}_k\in \Alg\bigl(\fP^{(n)}\bigr)$ if $A\in\Alg(\fP)$. 
It follows from $\widetilde{A}_k \eK\subseteq \eK$, for all $A\in \Alg(\fP)$ and all integers $1\leq k\leq n$,
that $\eK$ is a direct sum of subspaces from $\fP$. To see that $\eK$ is a direct sum of $n$ copies of the
same space, let $ 1\leq k, l\leq n$ be distinct integers and let operator $S_{k,l}\in \pB\bigl(\eX^{(n)}\bigr)$ 
be defined by $S_{k,l}(x_j)_{j=1}^{n}=(y_j)_{j=1}^{n}$, where $y_k=x_l$, $y_l=x_k$, and $y_j=x_j$ if 
$k\ne j\ne l$. Again, it is not hard to see that these operators are in $\Alg\bigl(\fP^{(n)}\bigr)$. It follows 
from $S_{k,l}\eK\subseteq \eK$ ($1\leq k\ne l\leq n$) that all direct summands of $\eK$ 
are the same subspace in $\eP$, that is, $\eK\in \eP^{(n)}$.
\end{proof}

Let us mention two things. First, by \Cref{cor02}, there exist reflexive pentagon subspace lattices of an arbitrary 
gap-dimension. Second, $\Lat\bigl(\Rk1(\Alg(\fP))\bigr)=\fP$ if $\fP$ is reflexive, however, $\fP$ is
not strongly reflexive. Hence, for a subspace lattice $\fL$, strong reflexivity is not equivalent to the equality
$\Lat\bigl(\Rk1(\Alg(\fL))\bigr)=\fL$.

The first proof of reflexivity of a pentagon subspace lattice with the gap-dimension $1$ was given by 
Halmos \cite{Hal}. In that paper, he asked whether every pentagon subspace lattice is reflexive. 
The answer was shown to be negative in \cite{LR} where a non-reflexive pentagon subspace lattice of gap-dimension $2$ is exhibited. More examples
of non-reflexive pentagon subspace lattices were given in \cite{KLL, LL}. Hence, although all
pentagon subspace lattices are isomorphic as lattices not all of them are reflexive, that is, reflexivity
is not a lattice theoretic property (on the other hand, strong reflexivity is).

\subsection{Double triangle}
The Hasse diagram of this abstract lattice $T$ is in \Cref{fig06}.

\begin{center}
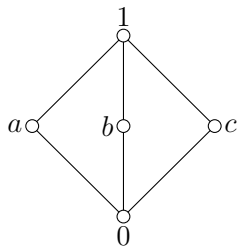

\begin{tikzpicture}[scale=1.2] 
\draw (0cm,0cm)--(1cm,1cm)--(0cm,2cm)--(-1cm,1cm)--(0cm,0cm)--(0cm,2cm);
\fill[white] (0cm,0 cm) circle (2pt);
\draw (0cm,0 cm) circle (2pt) node[below, scale=0.9] {$0$};
\fill[white] (-1cm,1cm) circle (2pt);
\draw (-1cm,1cm) circle (2pt) node[left, scale=0.9] {$a$};
\fill[white] (0cm,1cm) circle (2pt);
\draw (0cm,1cm) circle (2pt) node[left, scale=0.9] {$b$};
\fill[white] (1cm,1cm) circle (2pt);
\draw (1cm,1cm) circle (2pt) node[right, scale=0.9] {$c$};
\fill[white] (0cm,2cm) circle (2pt);
\draw (0cm,2 cm) circle (2pt) node[above, scale=0.9] {$1$};
\end{tikzpicture}
\captionof{figure}{The double triangle.}   \label{fig06}
\end{center}
Note that $T$ is a complemented lattice, however, $a, b$, and $c$ do not have unique complements.
If $\fT=\bigl\{ \{ 0\}, \eK, \eL,\eM, \eX\bigr\}$ is a realization of $T$ as a subspace lattice, then
$$ \eK\wedge\eL=\eK\wedge\eM=\eL\wedge\eM=\{ 0\}\quad \text{and}\quad
\eK\vee\eL=\eK\vee\eM=\eL\vee\eM=\eX. $$
We see that any two subspaces among $\eK$, $\eL$, and $\eM$ have a common complement. Hence,
$\dim(\eK)=\dim(\eL)=\dim(\eM)$ if $\eX$ is finite-dimensional. It follows that $T$ is realizable only in those
finite-dimensional Banach spaces which are of dimension $2n$ for a positive integer $n$.

We will show how one can realize $T$ as a subspace lattice.
Let $\eX$ be either a finite-dimensional Banach space of dimension $2n$ or an infinite-dimensional Banach space.
Assume that there exists a subspace $\eM_0\subseteq \eX$ that is algebraically complemented and is 
isomorphic to its complement, that is, there exists a subspace $\eM_\infty\subseteq \eX$ such that 
$\eM_0\cap \eM_\infty=\{ 0\}$, $\eM_0+\eM_\infty=\eX$ (i.e., $\eX=\eM_0\oplus \eM_\infty$),
and there exists an invertible operator $U\colon \eM_0\to \eM_\infty$.
If $\dim(\eX)=2n$, then $\eM_0$ and $\eM_\infty$ can be arbitrary $n$-dimensional subspaces of $\eX$ 
such that $\eM_0\cap \eM_\infty=\{ 0\}$.

Every vector $x\in \eX$ has a unique representation $x=x_0+x_\infty$, where $x_0\in \eM_0$ and 
$x_\infty\in \eM_\infty$. It follows that $Wx=U^{-1}x_\infty+Ux_0$ defines an operator on $\eX$. 
Consider vectors $x=x_0+x_\infty$ such that $Wx=x$. It is not hard to see that $x$ satisfies this condition 
if and only if $x_\infty=Ux_0$. It follows that $\Ker(I-W)=\{ x_0+Ux_0;\; x_0\in \eM_0\}$. Since $I-W$ is 
a bounded operator on $\eX$ its kernel $\eM_1=\Ker(I-W)$ is closed. It is not hard to see that 
$\eM_1$ is an algebraic complement of $\eM_0$.
Indeed, assume that $x\in \eM_0\cap\eM_1$. Then there exists $x_0\in \eM_0$ such that $x=x_0+Ux_0$.
Since $x-x_0\in \eM_0$ and $Ux_0\in \eM_\infty$ it follows from $\eM_0\cap\eM_\infty=\{ 0\}$ that $x=0$. 
To see that $\eM_0+\eM_1=\eX$, note that an arbitrary $x\in \eX$ is of the form
$x=x_0+x_\infty$ with $x_0\in \eM_0$ and $x_\infty\in \eM_\infty$. Since $U$ is a bijection, there exists 
a unique vector $z_0\in \eM_0$ such that $x_\infty=Uz_0$. It follows that $x=x_0+Uz_0=(x_0-z_0)+(z_0+Uz_0)$,
where $x_0-z_0\in \eM_0$ and $z_0+Uz_0\in \eM_1$.
A similar reasoning would show that $\eM_\infty\cap\eM_1=\{ 0\}$ and $\eM_\infty+\eM_1=\eX$.
Thus we may conclude that $\fT_U=\bigl\{ \{0\},\eM_0,\eM_1,\eM_\infty,\eX\bigr\}$ is a double triangle 
subspace lattice.

If $\fT=\bigl\{ \{0\},\eK,\eL,\eM,\eX\bigr\}$ is a double triangle subspace lattice in a finite-dimensional 
Banach space $\eX$ with $\dim(\eX)=2n$, then $\eK, \eL$, and $\eM$ are $n$-dimensional subspaces.
We have $\eX=\eK\oplus\eL$ and therefore every $x\in \eM$ has a unique representation $x=u+v$,
where $u\in \eK$ and $v\in \eL$. It is clear that $U\colon u\mapsto v$ defines an invertible operator 
$U\colon\eK\to\eL$. Hence, $\eM=\{ u+Uu;\; u\in \eK\}$. We conclude that $\fT$ is of the form described above,
that is, $\fT=\fT_U$.

Since $(\eM_0)_{-}=(\eM_1)_{-}=(\eM_\infty)_{-}=\eX$ the algebra $\Alg(\fT_U)$ does not contain 
rank-one operators, by \Cref{lem01}. The following lemma characterizes operators which are in $\Alg(\fT_U)$.
Note that we consider $\eX$ as the direct sum $\eM_0\oplus\eM_\infty$ and, for $x\in \eX$, the decomposition
$x=x_0+x_\infty$ means that $x_0\in \eM_0$ and $x_\infty\in \eM_\infty$.

\begin{lemma} \label{lem10}
An operator $S\in \pB(\eX)$ belongs to $\Alg(\fT_U)$ if and only if there exists an operator $A\in \pB(\eM_0)$
such that $Sx=S(x_0+x_\infty)=Ax_0+UAU^{-1}x_\infty$, for all $x\in \eX$.
\end{lemma}

\begin{proof}
If $A\in \pB(\eM_0)$ and we define $S\in \pB(\eX)$ by $Sx=S(x_0+x_\infty)=Ax_0+UAU^{-1}x_\infty$, 
for all $x\in \eX$, then $Sx_0=Ax_0\in \eM_0$ and $Sx_\infty=UAU^{-1}x_\infty\in \eM_{\infty}$, 
for all $x_0\in \eM_0$ and $x_\infty\in \eM_\infty$. Hence, $\eM_0$ and $\eM_\infty$ are invariant for $S$.
To see that $\eM_1$ is invariant for $S$, let $x=x_0+U x_0\in \eM_1$ be arbitrary. Then $Sx=Ax_0+UAx_0$. Since 
$Ax_0\in \eM_0$ we conclude that $Sx\in \eM_1$.

To prove the opposite implication, assume that $S\in \Alg(\fT_U)$. Since $\eM_0$ and $\eM_\infty$ are invariant for
$S$ and $\eX=\eM_0\oplus\eM_\infty$ we see that $S=A\oplus B$, where $A\in \pB(\eM_0)$ and 
$B\in \pB(\eM_\infty)$. Let $x_0\in \eM_0$ be arbitrary. Then $S(x_0+Ux_0)=Ax_0+BUx_0\in \eM_1$ as
$\eM_1$ is invariant for $S$. Hence, $BUx_0=UAx_0$, for all $x_0\in \eM_0$, and 
therefore $B=UAU^{-1}$.
\end{proof}

For every number $\omega\in \bC$, let $\eM_{\omega}=\{ x_0+\omega Ux_0;\; x_0\in \eM_0\}$. 
A similar reasoning as before for $\eM_1$ shows that $\eM_{\omega}$ are subspaces of $\eX$ such  that 
$\eX=\eM_{\omega_1}\oplus\eM_{\omega_2}$ if $\omega_1\ne \omega_2$
(this holds for $\omega_2=\infty$, as well).

\begin{lemma} \label{lem11}
If $\eN\in \Lat\Alg(\fT_U)$ is such that $\eN\cap \eM_\omega \ne \{ 0\}$, where $\omega\in \bC\cup\{\infty\}$,
then $\eM_\omega\subseteq \eN$. Hence, $\eN=\eX$ if there exist distinct $\omega_1$ and $\omega_2$
such that $\eN\cap \eM_{\omega_1} \ne \{ 0\}$ and $\eN\cap \eM_{\omega_2} \ne \{ 0\}$.
\end{lemma}

\begin{proof}
If $\eN\cap \eM_\omega\ne \{ 0\}$, for $\omega\ne \infty$, then there exists $0\ne x_0\in \eM_0$ such that
$x_0+\omega Ux_0\in \eN$. It follows that 
$(A\oplus UAU^{-1})(x_0+\omega Ux_0)=Ax_0+\omega UAx_0\in \eN$, for all $A\in \pB(\eM_0)$.
Since $\eM_0=\{ Ax_0;\; A\in \pB(\eM_0)\}$ we see that $\eM_\omega\subseteq \eN$. Similar reasoning
gives $\eM_\infty\subseteq \eN$ if $\eN\cap \eM_\infty\ne \{ 0\}$.
\end{proof}

It was proved by Halmos \cite{Hal} that a double triangle subspace lattice in a finite-dimensional subspace
is not reflexive. The following theorem slightly generalizes Halmos' result.

\begin{theorem} \label{theo04}
Subspace lattice $\fT_U$ is not reflexive since
$$\Lat\Alg(\fT_U)=\bigl\{ \eM_\omega;\; \omega\in \bC\cup\{\infty\} \bigr\}\cup \bigl\{ \{0\},\eX\bigr\}.$$
In particular, if $\fT$ is a double triangle subspace lattice in a finite-dimensional Banach space, then it is not reflexive.
\end{theorem}

\begin{proof}
Every $\eM_{\omega}$ is invariant for $\Alg(\fT_U)$. For $\omega\in \{ 0,1,\infty\}$, this is obvious since
$\fT_U\subseteq \Lat\Alg(\fT_U)$. For other $\omega$'s, let $x_0+\omega U x_0\in \eM_{\omega}$ be
arbitrary. Then $(A\oplus UAU^{-1})(x_0+\omega Ux_0)=Ax_0+\omega UAx_0\in \eM_{\omega}$,
for every $A\oplus UAU^{-1}\in \Alg(\fT_U)$. We have proved that
$\Lat\Alg(\fT_U)\supseteq \bigl\{ \eM_\omega;\; \omega\in \bC\cup\{\infty\} \bigr\}\cup \bigl\{ \{0\},\eX\bigr\}$.
To prove the opposite inclusion, let $\eN\in \Lat\Alg(\fT_U)$ be arbitrary. We may assume that $\eN\ne \{ 0\}$.

If $\eN\subseteq \eM_0$, then $\eN=\eM_0$, by \Cref{lem11}. Similarly, $\eN=\eM_\infty$, if 
$\eN\subseteq \eM_\infty$. Hence, we may suppose that $\eN\not\subseteq \eM_0$ and 
$\eN\not\subseteq \eM_\infty$. Let $x\in \eN\setminus \eM_0$ and $y\in \eN\setminus \eM_\infty$.
First, we consider the case that $x\in \eM_\infty$. In this case, $\eM_\infty\subseteq \eN$, by \Cref{lem11}.
We can write $y=y_0+y_\infty$, where $y_0\in \eM_0$ and $y_\infty\in \eM_\infty$. Hence, 
$0\ne y_0=y-y_\infty\in \eN$ and therefore $\eM_0\subseteq \eN$, by \Cref{lem11}. Thus, $\eN=\eX$,
in this case.

Suppose now that $x\not\in \eM_\infty$. Then $x=x_0+x_\infty$, where $0\ne x_0\in \eM_0$ and
$0\ne x_\infty\in \eM_\infty$. Let $0\ne z_0=U^{-1}x_\infty\in \eM_0$. If $x_0$ and $z_0$ are linearly
independent, then for every pair of vectors $u_0\in \eM_0$ and $u_\infty\in \eM_\infty$, there exists
an operator $A\in \pB(\eM_0)$ such that $Ax_0=u_0$ and $UAz_0=u_\infty$. It follows that
$u_0+u_\infty=(A\oplus UAU^{-1})(x_0+x_\infty)\in \eN$, that is, $\eN=\eX$. It remains to consider the case
when $x_0$ and $z_0$ are linearly dependent, say $z_0=\omega x_0$. Note that $\omega\ne 0$ since
$x_0\ne 0$ and $z_0\ne 0$. Thus, $x=x_0+\omega Ux_0\in \eM_\omega$. By \Cref{lem11},
$\eM_\omega\subseteq \eN$. If $\eM_\omega\ne \eN$, then there exists $w\in \eN\setminus \eM_\omega$.
Since $\eX=\eM_0+\eM_\omega$ we can write $w=s_0+(t_0+\omega Ut_0)$, where $0\ne s_0\in \eM_0$
and $t_0+\omega Ut_0\in \eM_\omega $ (i.e., $t_0\in \eM_0$). Vectors $w$ and $t_0+\omega Ut_0$
are in $\eN$, hence, $s_0\in \eN$. By \Cref{lem11}, $\eM_0\subseteq \eN$. We conclude that $\eN=\eX$.
\end{proof}

 It is not known if there exists a reflexive double triangle subspace lattice. If 
$\fT=\bigl\{ \{0\},\eK,\eL,\eM,\eX\bigr\}$ is a reflexive double triangle subspace lattice in $\eX$, then,
of course, $\eX$ must be infinite-dimensional and none of the sums $\eK+\eL$, $\eK+\eM$, $\eL+\eM$
is closed (see \cite[Theorem 1]{Lon2} for the proof in the Hilbert space case --- a simple modification gives
a proof for a general Banach space). Kissin \cite{Kis} showed that $\fT$ (in a complex Hilbert space) is non-reflexive 
if $\Alg(\fT)$ contains a non-zero finite-rank operator.

\subsection*{Acknowledgment}
The author is grateful to Professors Bojan Kuzma and Marko Kandi\'{c} for reading the first version of this
manuscript and for their valuable comments. The work on this survey paper was partially supported by the 
Slovenian Research and Innovation Agency through the research program P2-0268.

	

\begin{thebibliography}{99}

\bibitem{Arv} 
	W. Arveson, 
	\textit{Operator algebras and invariant subspaces},
	Ann. of Math. (2) 100 (1974), 433--532.
	
\bibitem{Con} 
	J. B. Conway, 
	\textit{A complete Boolean algebra of subspaces which is not reflexive},
	Bull. Amer. Math. Soc. 79 (1973), 720--722.

\bibitem{Dav} 
	K. R. Davidson, 
	\textit{Commutative subspace lattices},
	Indiana Univ. Math. J. 27 (1978), 479--490.

\bibitem{DWY} 
	A. Dong, W. Wu, W. Yuan, 
	\textit{On small subspace lattices in Hilbert space},
	J. Austral. Math. Soc. Ser. 96 (2014), 44--60.
		
\bibitem{Erd} 
	J. A. Erdos, 
	\textit{Operators of finite rank in nest algebras},
	J. London Math. Soc. 43 (1968), 391--397.
		
\bibitem{HLR} 
	D. W. Hadwin, W. E. Longstaff, P. Rosenthal,
	\textit{Small transitive lattices},
	Proc. Amer. Math. Soc. 87 (1983), 121--124.
			
\bibitem{Hal1} 
	P. R. Halmos, 
	\textit{Ten problems in Hilbert space},
	Bull. Amer. Math. Soc. 76 (1970), 887--933.

\bibitem{Hal} 
	P. R. Halmos, 
	\textit{Reflexive lattices of subspaces},
	J. London Math. Soc. (2) 4 (1971), 257--263.

\bibitem{Hal2} 
	P. R. Halmos, 
	\textit{Ten years in Hilbert space},
	Integral Equations Operator Theory 2 (1979), 529--564.
	
\bibitem{Har} 
	K. J.  Harrison, 
	\textit{Certain distributive lattices of subspaces are reflexive},
	J. London Math. Soc. (2) 8 (1974), 51--56.

\bibitem{HL} 
	K. J.  Harrison, W. E. Longstaff,
	\textit{Reflexive subspace lattices in finite-dimensional Hilbert spaces},
	Indiana Univ. Math. J. 26 (1977), 1019--1025.

\bibitem{Hay} 
	R. Haydon, 
	\textit{Reflexivity of commutative subspace lattices},
	Proc. Amer. Math. Soc. 115 (1992), 1057--1060.
	
\bibitem{HM} 
	A. Hopenwasser, R. Moore 
	\textit{Finite rank operators in reflexive operator algebras},
	J. London Math. Soc. (2) 27 (1983), 331--338.
		
\bibitem{Joh} 
	R. E. Johnson,
	\textit{Distinguished rings of linear transformations},
	Trans. Amer. Math. Soc. 111 (1964), 400--412.
	
\bibitem{KLL} 
	A. Katavolos, M. S. Lambrou, W. E. Longstaff,
	\textit{Pentagon subspace lattices on Banach spaces},
	J. Operator Theory 46 (2001), 355--380.

\bibitem{Kis} 
	E. Kissin,
	\textit{On some reflexive lattices of subspaces},
	J. Operator Theory 25 (1991), 141--162.
	
\bibitem{Lam} 
	M. S. Lambrou, 
	\textit{Approximants, commutants and double commutants in normed algebras},
	J. London Math. Soc. (2) 25 (1982), 499--512.

\bibitem{Lam1} 
	M. S. Lambrou,
	\textit{On the rank of operators in reflexive algebras},
	Linear Algebra Appl. 142 (1990), 211--235.
		
\bibitem{LL} 
	M. S. Lambrou, W. E. Longstaff,
	\textit{Non-reflexive pentagon subspace lattices},
	Studia Math. (2) 125 (1997), 187--199.

\bibitem{LauL} 
	C. Laurie, W. E. Longstaff,
	\textit{A note on rank-one operators in reflexive algebras},
	Proc. Amer. Math. Soc. 89 (1983), 293--297.
	
\bibitem{LT} 
	J. Lindenstrauss, L. Tzafriri,
	\textit{On the complemented subspace problem},
	Israel J. Math. 9 (1971), 263--269.
	
\bibitem{Lon1} 
	W. E. Longstaff, 
	\textit{Strongly reflexive lattices},
	J. London Math. Soc. (2) 11 (1975), 491--498.

\bibitem{Lon2} 
	W. E. Longstaff, 
	\textit{Non-reflexive double triangles},
	J. Austral. Math. Soc. Ser. A35 (1983), 349--356.

\bibitem{Lon3} 
	W. E. Longstaff, 
	\textit{Operators of rank one reflexive algebras},
	Canadian J. Math. 28 (1976), 19--23.

\bibitem{LR} 
	W. E. Longstaff, P. Rosenthal,
	\textit{On two questions of Halmos concerning subspace lattices},
	Proc. Amer. Math. Soc. 75 (1979), 85--86.

\bibitem{Mac} 
	G. W. Mackey,
	\textit{Note on a theorem of Murray},
	Bull. Amer. Math. Soc. 52 (1946), 322--325.

\bibitem{Mur} 
	F. J. Murray,
	\textit{Quasi-complements and closed projections in reflexive Banach spaces},
	Trans. Amer. Math. Soc. 58 (1945), 77--95.

\bibitem{Rin1} 
	J. R. Ringrose, 
	\textit{On some algebras of operators},
	Proc. London Math. Soc. 15 (1965), 61--83.

\bibitem{Rin2} 
	J. R. Ringrose, 
	\textit{On some algebras of operators. II},
	Proc. London Math. Soc. 16 (1966), 385--402.
	
\bibitem{Shu} 
	V. S. Shulman, 
	\textit{Projection lattices in a Hilbert space},
	Funct. Anal. Appl. 23 (1989), 158--159.

\end{thebibliography}
\end{document}